\documentclass[preprint,12pt]{elsarticle}

\PassOptionsToPackage{table}{xcolor}
\PassOptionsToPackage{colorlinks}{hyperref}

\usepackage{amssymb}
\usepackage{amsmath}

\usepackage{graphicx}
\usepackage{svg}
\usepackage{subcaption}
\usepackage{tabularx}
\usepackage{booktabs}
\usepackage{enumitem}
\usepackage{xurl}

\usepackage{xcolor}
\usepackage{colortbl}

\usepackage{algorithm}
\usepackage{algpseudocode}

\usepackage{physics}

\usepackage{hyperref}
\usepackage{cleveref}

\crefname{assumption}{Assumption}{Assumptions}
\crefname{definition}{Definition}{Definitions}
\crefname{lemma}{Lemma}{Lemmas}
\crefname{remark}{Remark}{Remarks}
\crefname{theorem}{Theorem}{Theorems}
\crefname{proposition}{Proposition}{Propositions}
\crefname{section}{Section}{Sections}
\crefname{figure}{Figure}{Figures}
\crefname{equation}{}{}
\crefname{table}{Table}{Tables}
\crefname{appendix}{}{}
\crefname{algorithm}{Algorithm}{Algorithms}

\journal{Journal of Computational and Applied Mathematics}

\begin{document}

  \begin{frontmatter}

  \title{High-Order Symmetric Positive Interior Quadrature Rules on Two and Three Dimensional Domains}

  \author[uw]{Moustapha Diallo\corref{cor1}} 
  \ead{mdiallo@uwaterloo.ca}
  \author[uw]{Zelalem Arega Worku}
  \ead{zelalem.worku@uwaterloo.ca}
  \cortext[cor1]{Corresponding author.} 
  \affiliation[uw]{organization={University of Waterloo, Department of Applied Mathematics},
              addressline={200 University Avenue West}, 
              city={Waterloo},
              postcode={N2L 3G1}, 
              state={ON},
              country={Canada}}

  \begin{abstract}
  Fully symmetric positive interior (f‑SPI) quadrature rules are key building blocks for high‑order discretizations of partial differential equations, yet high‑degree rules with few nodes remain scarce on reference elements commonly used in mesh generation. We construct new f‑SPI rules on the square, cube, prism, and pyramid by coupling a variable parameterization that enforces positivity and interiority with an efficient Levenberg–Marquardt optimization and a symmetry‑aware node‑reduction strategy that eliminates and collapses orbits, allowing transitions between symmetry types. The resulting rules achieve degrees up to 77 on the square, 45 on the cube, and 30 on the prism and pyramid, and for most degrees use fewer nodes than previously published f‑SPI quadrature rules. Verification tests demonstrate comparable accuracy to existing rules. Complete node and weight data are also provided.
  \end{abstract}


  \end{frontmatter}

  \section{Introduction}
  \label{intro}
  Quadrature rules are approximation tools for definite integrals. They play an important role in the process of solving partial differential equations (PDEs) numerically. For example, they are heavily used in discontinuous Galerkin (DG) finite element methods \cite{osti_4491151,COCKBURN1998199,hesthaven2008} and summation-by-parts (SBP) discretizations \cite{fernandez2014review,svard2014review,hicken2016multidimensional}. The mathematical properties of such numerical methods, e.g., accuracy and stability, are affected by the type of quadrature rule used. In particular, high-degree fully symmetric, positive, and interior (f-SPI) rules with minimal nodes are often sought to obtain more accurate, stable, and efficient (speed and memory) numerical methods. The degree of available quadrature rules can also be a limiting factor on the order of numerical methods. In the case of DG methods, a degree $q$ method can require a quadrature rule of degree up to $6q$ \cite{KIRBY2003249,PERSSON20091585,Williams2019EntropyL2DG}. The construction of very high order f-SPI rules is thus a critical endeavor. Triangles, quadrilaterals, tetrahedra, pyramids, prisms, and hexahedra often arise as discretization elements in the application of numerical methods. For triangles and tetrahedra, there exist very high order f-SPI quadrature rules with relatively few nodes, e.g., see \cite{Worku2026,chuluunbaatar2022,Jaskowiec2021main}; hence, our work focuses on the other four domains.

  The construction of quadrature rules consists of two main steps. The first step is initialization of the rule where a distribution of nodes and weights is chosen. Several strategies have been explored in previous work, including simulated \cite{WANDZURAT20031829}, particle swarm optimization \cite{worku2024quadrature}, random distribution \cite{WITHERDEN20151232,Jaskowiec2021main,jaskowiec2021addendum}, tensor-product \cite{SLOBODKINS2023229,XIAO2010663,Worku2026}, and the use of Padua points \cite{FESTA20124296}. The second step is refinement of the rule where a solver is used to find optimal distributions of nodes and weights. Commonly used techniques include Newton's method \cite{XIAO2010663,FESTA20124296,Jaskowiec2021main,jaskowiec2021addendum,SLOBODKINS2023229,WANDZURAT20031829} and the Levenberg-Marquardt algorithm (LMA) \cite{Worku2026,WITHERDEN20151232,kubatko2013}, often with some modifications to enforce domain constraints. Nodes are also added or removed at this stage. With random initializations, the intent usually is to start with a minimal node quadrature rule, but this is sometimes too ambitious and additional nodes are needed to satisfy the moment equations. In contrast, tensor-product initializations, which can easily satisfy the moment equations, lead to more nodes than necessary, requiring the removal of a relatively large number of nodes to reduce node count. The best results (in terms of degree and node count) we found on the cube, prism, and pyramid go up to degree 21, 20, and 20, respectively \cite{Jaskowiec2021main,jaskowiec2021addendum}. They were constructed using node addition followed by node elimination. Similarly, the best results we can find for the square go up to degree 21 \cite{WITHERDEN20151232}.

  Our work builds on the methods and results from \cite{Worku2026}, where a combination of tensor-product-based initialization and node elimination is used to build very high order quadrature rules on the triangle (up to degree 84) and tetrahedron (up to degree 40). In this work, we introduce two new ways of initializing quadrature rules on the pyramid with relatively few nodes, one of which uses results from \cite{Worku2026}. We also present an improvement on the node elimination approach for symmetric rules, and implement a variable parameterization that turns the quadrature moment problem into an unconstrained optimization problem, enabling the use of the LMA without constraints. As a result, we push the attainable degrees upward while reducing node counts relative to existing results, obtaining rules up to degree 77, 45, 30, and 30 for the square (quadrilateral), cube (hexahedron), prism, and pyramid, respectively. 

  The paper is structured as follows. \cref{problem} defines the problem of finding a degree $q$ quadrature rule, and \cref{lma} presents the solver used. \cref{initiali-quadratures,node-elimination} explain the construction of initial quadrature rules and the node elimination process, respectively. Finally, \cref{results} details the results, which are followed by the conclusions in \cref{conclusions}.

  \section{Problem}
  \label{problem}

\begin{figure}[!t]
      \centering

      \begin{subfigure}{0.325\linewidth}
        \centering
        \includegraphics[
          width=\linewidth
        ]{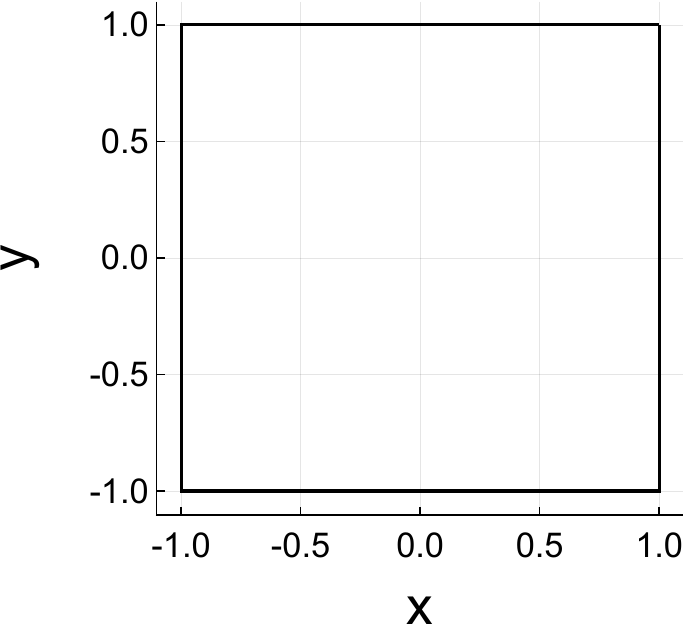}
        \caption{Square}
        \label{fig:square}
      \end{subfigure}\hfill
      \begin{subfigure}{0.45\linewidth}
        \centering
        \includegraphics[
          width=\linewidth,
        ]{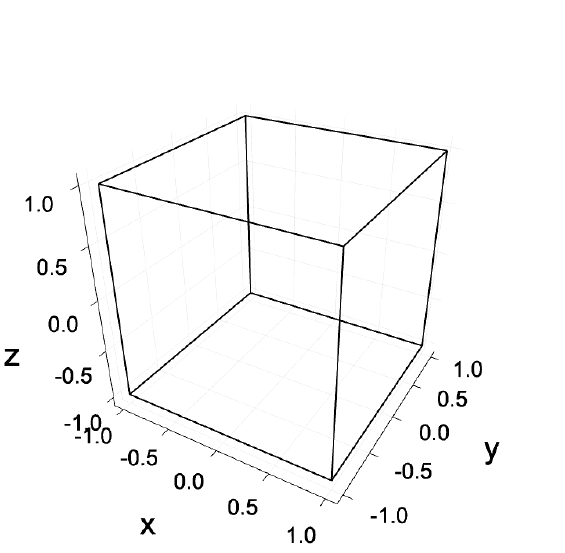}
        \caption{Cube}
        \label{fig:cube}
      \end{subfigure}
      \begin{subfigure}{0.50\linewidth}
        \centering
        \includegraphics[
          width=\linewidth,
          ]{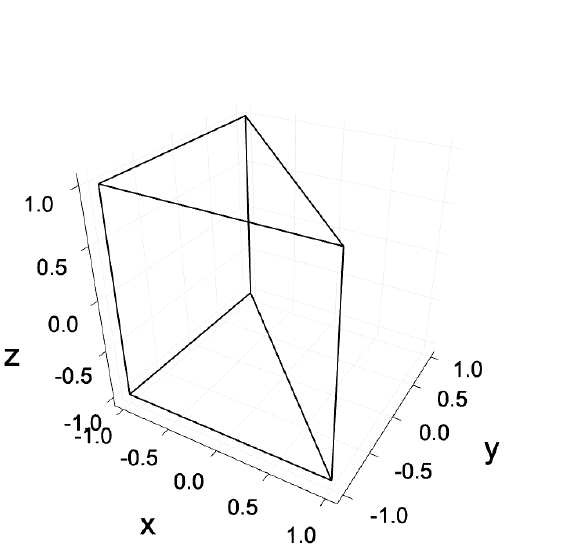}
        \caption{Prism}
        \label{fig:prism}
      \end{subfigure}\hfill
      \begin{subfigure}{0.50\linewidth}
        \centering
        \includegraphics[
          width=\linewidth,
        ]{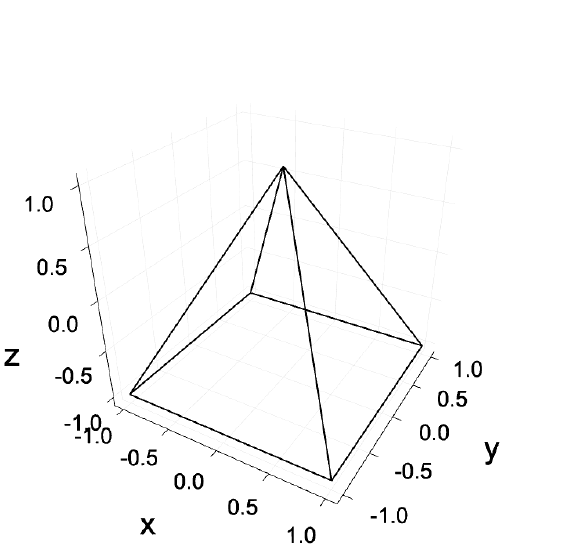}
        \caption{Pyramid}
        \label{fig:pyramid}
      \end{subfigure}

      \caption{Visual of the domains.}
      \label{fig:domains}
  \end{figure}

  An f-SPI quadrature rule of degree $q$ in a domain $D$ is a set of nodes and weights $Q = \{(\mathbf{x}_j, w_j) : \mathbf{x}_j \in D, \ w_j \in (0,\infty), \ j = 1,2,..,n \}$ used to approximate integrals over $D$. The nodes are fully symmetrically distributed within the domain, and the approximation is exact for all polynomials of degree at most $q$, $P_q$; that is, 
  \begin{equation}
      \label{eq:deg-q-rule}
      \int_D f(\mathbf{x})\dd\mathbf{x} = \sum_{j=1}^n w_jf(\mathbf{x}_j), \ \ \forall \ f \in P_q.
  \end{equation}
  Our domains of interest are defined as follows: 
  \begin{itemize}
      \item Square (quadrilateral): $sqr = {\{(x,y)\mid |x|,|y|<1\}},$
      \item Cube (hexahedron): $cube = {\{(x,y,z)\mid |x|,|y|,|z|<1\}},$
      \item Prism: $pri = {\{(x,y,z)\mid -1<y<-x<1, \ |z|<1\}},$
      \item Pyramid: $pyr = {\{(x,y,z)\mid |z|<1, \ |x|,|y|<(1-z)/2\}},$
  \end{itemize}
  which are also shown pictorially in \cref{fig:domains}.

  Suppose $P_q$ is spanned by the basis $B_q = \{p_i \mid p_i \in P_q, \ i = 1,2,\ldots,m\}$, then the problem of finding a degree $q$ quadrature rule over the domain $D$ can be formulated as the following minimization problem:
  \begin{equation}\label{eq:minimization_prob}
    \min_{\{\mathbf{x}_j,w_j\}_{j=1}^n} \|\mathbf{f}-V^{T}\mathbf{w}\|_2 \quad \text{s.t.}\quad \mathbf{x}_j\in D,\quad w_j>0,\quad j=1,\ldots,n,
  \end{equation}
  where $f_i = \int_D p_i(\mathbf{x})\dd\mathbf{x}$ and $V_{ji} = p_i(\mathbf{x}_j)$ for all $i = 1,\ldots,m$ and $j = 1,\ldots,n$. We use an orthogonal polynomial basis to obtain a well-conditioned Vandermonde matrix, $V$, and an easy to compute vector of the exact polynomial integrals (${f_i = \delta_{i1}*\mathrm{constant}}$). In theory, the optimality condition of the problem is that the objective is zero, but we terminate when $\|\mathbf{f}-V^{T}\mathbf{w}\|_2 < 10^{-66}$. 

  There are computational advantages to solving \cref{eq:minimization_prob} for fully symmetric quadrature rules. First, the number of degrees of freedom can be significantly smaller than the number of node coordinates. A symmetric distribution leads to dependent nodes of the same weight. Such a collection of nodes can be compressed into an orbit parameterized by a representative node and weight. Orbits are classified into different symmetries based on the number and relationship of nodes. The number of nodes associated with each orbit in a symmetry represents the size of that symmetry. Table~\ref{tab:symmetries} displays pairs of (symmetry, symmetry size) for each domain, grouped by the number of parameters in the symmetries. Second, symmetry implies that many moment conditions are satisfied automatically. For example, symmetric quadratures on the square exactly integrate any polynomial that is odd in at least one of the variables. Those polynomials can be removed from the orthogonal basis to reduce the number of columns of the Vandermonde matrix. The resulting set of polynomials is called an objective basis. More details on the symmetries and the objective bases used can be found in \cite{WITHERDEN20151232}.

  \begin{table}[!tbp]
      \caption{Pairs of (symmetry, symmetry size) grouped by the number of parameters in the symmetries.}
      \label{tab:symmetries}
      \centering
      \small
      \begin{tabular}{l|l|l|l|l}
          & 0-param & 1-param & 2-param & 3-param \\ \hline
          Square & $(S_1,1)$ & $(S_2,4) \ (S_3,4)$ & $(S_4,8)$ & \\ \hline
          Cube & $(S_1,1)$ & $(S_2,6), \ (S_3,8), \ (S_4,12)$ & $(S_5,24), \ (S_6,24)$ & $(S_7,48)$ \\ \hline
          Prism & $(S_1,1)$ & $(S_2,2), \ (S_3,3)$ & $(S_4,6), \ (S_5,6)$ & $(S_6,12)$ \\ \hline
          Pyramid &  & $(S_1,1)$ & $(S_2,4), \ (S_3,4)$ & $(S_4,8)$ \\
      \end{tabular}
      
  \end{table}

  \section{Levenberg-Marquardt Algorithm}
  \label{lma}

  We use the LMA to solve the minimization problem stated in  \cref{eq:minimization_prob}. Similar to Newton's method, the LMA starts with an initial state (of the variables) $\mathbf{v}^{(0)}$ and generates a sequence of states, $\mathbf{v}^{(k)}$, 
  \begin{equation}
      \mathbf{v}^{(k+1)} = \mathbf{v}^{(k)} + \mathbf{dv}^{(k)},
  \end{equation}
  that ideally converges to a global minimum of the problem.
  There are several variants of the algorithm, and we use the one described in \cite{Worku2026}. The LMA is designed for unconstrained optimization problems, but \cref{eq:minimization_prob} has constraints --- interior nodes and positive weights. We explored two ways of addressing this issue.

  The first approach consists of using a Cartesian parameterization of the variables and projection of $\mathbf{v}^{(k+1)}$ inside the domain. We perform the projection through scaling down $\mathbf{dv}^{(k)}$ by the largest scaling factor $\alpha \in (0, 1]$ that ensures $\mathbf{v}^{(k+1)}$ lies within the domain. This scalar is computed as follows:
  \begin{align}
      \alpha &= \min{\{\alpha_i: i =1,\ldots,n_v\}}, 
      \nonumber
      \\
      \alpha_i &=
      \begin{cases}
          \frac{v_i^{(k)} + dv_i^{(k)} - (v_{i,\mathrm{min}} + \epsilon_m)}{dv^{(k)}_{i}}, & \text{if } v_i^{(k)} + dv_i^{(k)} \leq v_{i,\mathrm{min}} \\ 
          \frac{v_i^{(k)} + dv_i^{(k)} - (v_{i,\mathrm{max}} - \epsilon_m)}{dv^{(k)}_{i}}, & \text{if } v_i^{(k)} + dv_i^{(k)} \ge v_{i,\mathrm{max}} \\
          1, & \text{otherwise},
      \end{cases}
  \end{align}
  where $n_v$ is the number of variables, $(v_{i,\mathrm{min}},v_{i,\mathrm{max}})$ is the domain of variable $v_i$, and $\epsilon_m$ is the machine precision for the number type being used, \verb|Float64| or \verb|BigFloat|. Coordinate variables are normalized so that their domain is $(0,1)$, and weight variables have domain $(0,\infty)$. A similar projection method is used in \cite{Worku2026}.

  The second approach reformulates \cref{eq:minimization_prob} as an unconstrained optimization problem using exponential parameterization of the variables. This reparameterization is particularly useful because it enforces positivity and interiority by construction. For a variable with Cartesian parameter $x\in (0,1)$, we use the exponential parameter $z \in \mathbb{R}$ given by the following:
  \begin{equation}\label{eq:exp1}
      x = \frac{1}{1 + e^{-sz}}, \quad s = 0.01.
  \end{equation}
  For a variable with Cartesian parameter $x\in (0,\infty)$, we use the exponential parameter $z \in \mathbb{R}$ given by the following:
  \begin{equation}\label{eq:exp2}
      x = e^{sz}.
  \end{equation}
  We experimented with several values for $s$ and selected $s=0.01$, which consistently led to the elimination of the largest number of nodes.

  Exponential parameterization is computationally faster because it avoids repeated projections, but it is more prone to numerical instabilities, often resulting in ill-conditioned or singular matrices in the LMA update. To balance robustness and efficiency, we employ a hybrid strategy: exponential parameterization is used by default, and the solver reverts locally to Cartesian parameters whenever an inversion failure is detected. 

  \section{Initial Quadrature Construction}
  \label{initiali-quadratures}

  The first step of our approach is to construct the initial quadrature rules. At this stage, the goal is to make quadrature rules that have few nodes without being too restrictive, e.g. having no nodes on some symmetries. We mainly construct them by decomposing the domains into lower-dimensional components and using the best known quadratures on those components.

  Suppose that we are constructing a degree $q$ quadrature rule, and let $p_q(\mathbf{x})$ be an arbitrary polynomial of degree at most $q$. We present the approach for each domain considered.

  \subsection{Square} 
  The square can be viewed as the tensor-product of two line segments on $(-1,1)$:
  \begin{align*}
      \int_{sqr}p_q(\mathbf{x})d\mathbf{x} &= \int_{-1}^1\left(\int_{-1}^1p_q(\mathbf{x})dx\right)dy \\ 
      & = \int_{-1}^1\sum_{a=1}^A w_{la}p_q(x_{la},y)dy\\
      & = \sum_{a=1}^A \sum_{b=1}^A w_{la}w_{lb}p_q(x_{la},x_{lb}),
  \end{align*}
  where $\{(x_{la}, w_{la})\mid a=1,..,A\}$ is the Gauss-Legendre quadrature rule on $(-1,1)$ of degree at least $q$ with the minimal odd number of nodes. The resulting quadrature on the square is 
  $\{(x_{la},x_{lb},w_{la}w_{lb})\mid a,b=1,\ldots,A\}$. The choice of Gauss-Legendre quadratures with an odd number of nodes ensures that the initial quadrature has nodes on the $S_1$ and $S_2$ symmetries.

  \subsection{Cube} 
  The cube can be viewed as the tensor-product of three line segments on $(-1,1)$:
  \begin{align*} 
      \int_{cube}p_q(\mathbf{x})d\mathbf{x} & = \int_{-1}^1\left(\int_{-1}^1\left(\int_{-1}^1p_q(\mathbf{x})dx\right)dy\right)dz \\ 
      & = \int_{-1}^1\left(\int_{-1}^1\sum_{a=1}^A w_{la}p_q(x_{la},y,z)dy\right)dz\\
      & = \int_{-1}^1\sum_{a=1}^A \sum_{b=1}^A w_{la}w_{lb}p_q(x_{la},x_{lb},z)dz \\ 
      & = \sum_{a=1}^A\sum_{b=1}^A \sum_{c=1}^A w_{la}w_{lb}w_{lc}p_q(x_{la},x_{lb},x_{lc}).
  \end{align*}
  The resulting quadrature on the cube is $\{(x_{la},x_{lb},x_{lc}, w_{la}w_{lb}w_{lc})\mid a,b,c=1,\ldots,A\}$. Here, the use of an odd degree Legendre-Gauss quadrature rule on $(-1,1)$ ensures that the cube has nodes on its $S_1$, $S_2$, $S_4$, and $S_5$ symmetries.

  \subsection{Prism} 
  The prism can be viewed as the tensor-product of the line segment $(-1,1)$ and the triangle $tri=\{(x,y)\mid -1<y<-x<1\}$:
  \begin{align*}
      \int_{pri}p_q(\mathbf{x})d\mathbf{x} & = \int_{-1}^1\left(\int_{tri}p_q(\mathbf{x})dxdy\right)dz \\ & = \int_{-1}^1\sum_{b=1}^B w_{tb}p_q(x_{tb},y_{tb},z)dz \\ & =
      \sum_{b=1}^B\sum_{a=1}^A w_{tb}w_{la}p_q(x_{tb},y_{tb},x_{la}),
  \end{align*}
  where $\{(x_{tb},y_{tb}, w_{tb})\mid b=1,\ldots,B\}$ is the degree $q$ f-SPI quadrature rule presented in \cite{Worku2026} for the triangle. We add the centroid whenever absent to ensure that initial quadrature rules have nodes on the $S_1$ and $S_2$ symmetries. The resulting quadrature rule on the prism is $\{(x_{tr},y_{tr},x_{ls}, w_{tr}w_{ls})\mid {r=1,\ldots,b},\ {s=1,\ldots,a}\}$.

  \subsection{Pyramid} 
  We developed two ways of initializing the pyramid. One uses an algebraic construction while the other uses a geometric construction. The algebraic construction consists of expressing integration over the pyramid as a composition of an integral over $(-1,1)$ and an integral over $(-1,1)^2$:
  \begin{align*}
      \int_{pyr}p_q(\mathbf{x})d\mathbf{x} & = 
      \int_{-1}^1\left(\int_{-\alpha_z}^{\alpha_z}\int_{-\alpha_z}^{\alpha_z} p_q(x,y,z)dxdy\right)dz \\ & = \int_{-1}^1\left(\alpha_z^2\int_{-1}^1\int_{-1}^1p_q(\alpha_z x,\alpha_z y,z)dxdy\right)dz \\ & = \int_{-1}^1\alpha_z^2\sum_{c=1}^{C} w_{sc}p_q(\alpha_zx_{sc},\alpha_zy_{sc},z)dz \\ & = \sum_{e=1}^{E}\sum_{c=1}^C\alpha_{x_{le}}^2w_{sc}w_{le}p_q(\alpha_{x_{le}}x_{sc},\alpha_{x_{le}}y_{sc},x_{le}),
  \end{align*}
  where $\alpha_z = (1-z)/2$, $\{(x_{le},w_{le})\mid e=1,\ldots,E\}$ is the degree $q+2$ Gauss-Legendre quadrature rule on $(-1,1)$, and $\{(x_{sc},y_{sc},w_{sc})\mid c=1,\ldots,C\}$ is our best result for degree $q$ on the square with the center included to ensure that the initial quadrature on the pyramid has nodes on the $S_1$ symmetry. The resulting quadrature rule on the pyramid is $\{(\alpha_{x_{le}}x_{sc},\alpha_{x_{le}}y_{sc},x_{le},\alpha_{x_{le}}^2w_{sc}w_{le})\mid c=1,\ldots,C,\ e=1,\ldots,E\}$.

\begin{figure}[!t]
      \begin{subfigure}{0.24\linewidth}
        \centering
        \includegraphics[
          width=\linewidth, 
          ]{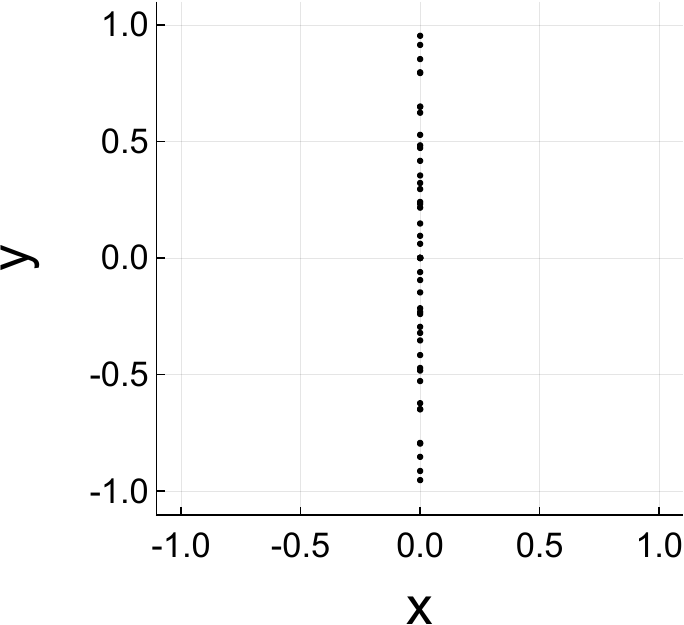}
        \caption{}
        \label{fig:pyr-geoA}
      \end{subfigure}\hfill
      \begin{subfigure}{0.24\linewidth}
        \centering
        \includegraphics[
          width=\linewidth, 
        ]{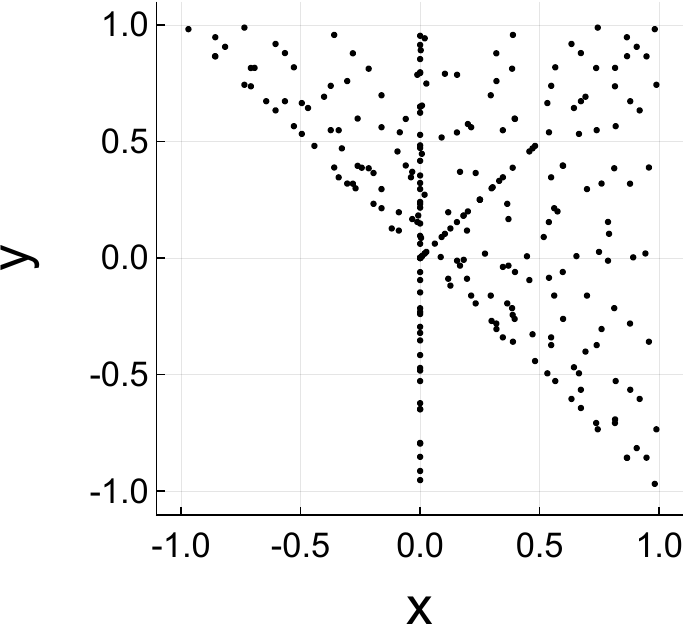}
        \caption{}
        \label{fig:pyr-geoB}
      \end{subfigure}\hfill
      \begin{subfigure}{0.24\linewidth}
        \centering
        \includegraphics[
          width=\linewidth, 
        ]{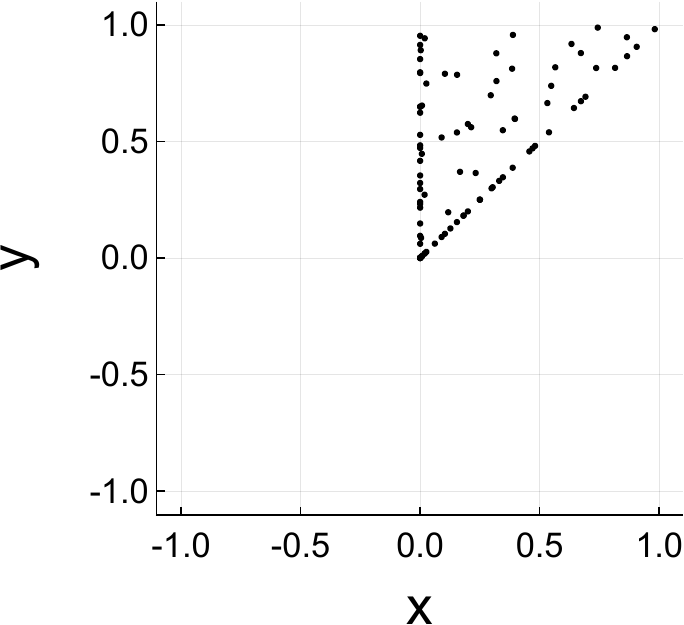}
        \caption{}
        \label{fig:pyr-geoC}
      \end{subfigure}\hfill
      \begin{subfigure}{0.24\linewidth}
        \centering
        \includegraphics[
          width=\linewidth, 
        ]{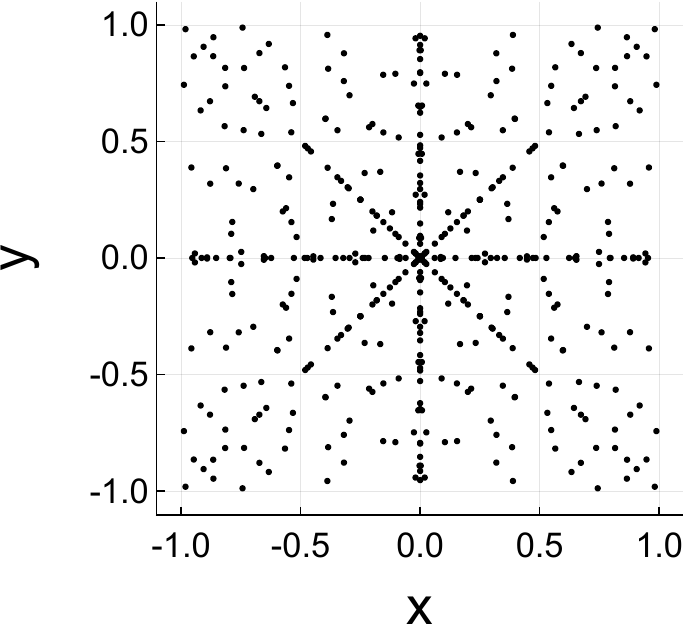}
        \caption{}
        \label{fig:pyr-geoD}
      \end{subfigure}\hfill

      \caption{Base view of the geometric construction of a degree 15 quadrature rule on the pyramid.}
      \label{fig:pyr-geo15construction}
  \end{figure}

  The geometric construction splits the pyramid into triangles and tetrahedra. The steps involved are:
  \begin{enumerate}[label=(\alph*)]
      \item map the nodes of the degree $q$ triangular quadrature rule from \cite{Worku2026} to $\{(0,y,z)\mid |z|<1, \ |y|<(1-z)/2\} \subset pyr$;
      \item map the nodes of the degree $q-2$ tetrahedral quadrature rule from \cite{Worku2026} to $\{(x,y,z)\mid |z|<1, \ |x|,|y|<(1-z)/2, \ -x<y\} \subset pyr$;
      \item remove the nodes not in $\{(x,y,z)\mid 0 \leq x \leq y < 1\, \ |z|<1\}$;
      \item use the symmetries of the pyramid to generate more nodes;
      \item set each weight to (pyramid volume)/(node count), then pass the quadrature rule to the LMA.
  \end{enumerate}
  \cref{fig:pyr-geo15construction} provides a visual of steps (a) to (d) for the geometric construction of a degree 15 quadrature rule on the pyramid. 

  It appears that the only previously used method to make a convergent quadrature rule, without using an optimization algorithm, on the pyramid is through a Duffy transformation \cite{duffy1982}. It works by mapping an integral over the pyramid to an integral over the cube but leads to a quadrature rule with a very high number of nodes as a degree $q$ rule on the pyramid requires a degree $2q + 2$ rule on the cube. The algebraic construction described above is also convergent without optimization; however, it only requires a degree $q$ rule on the square and a degree $q+2$ Gauss-Legendre rule, which results in a significantly lower node count compared to the Duffy construction. For instance, a degree 15 quadrature rule on the pyramid has 1787 and 441 nodes with the Duffy and algebraic constructions, respectively. The Duffy transformation in this case uses the degree 33 quadrature rule we have obtained on the cube to produce the quadrature rule on the pyramid with 1787 nodes.

\begin{figure}[!t]
      \begin{subfigure}{0.5\linewidth}
        \centering
        \includegraphics[
          width=\linewidth, 
        ]{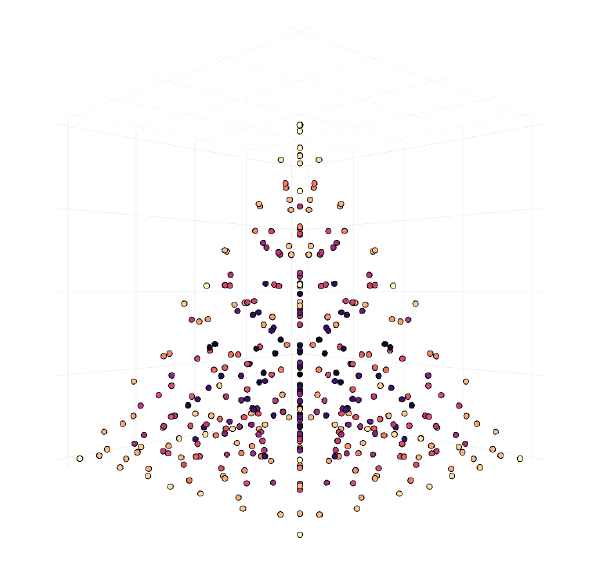}
        \caption{Geometric}
        \label{fig:pyr-geo15}
      \end{subfigure}\hfill
      \begin{subfigure}{0.5\linewidth}
        \centering
        \includegraphics[
          width=\linewidth, 
        ]{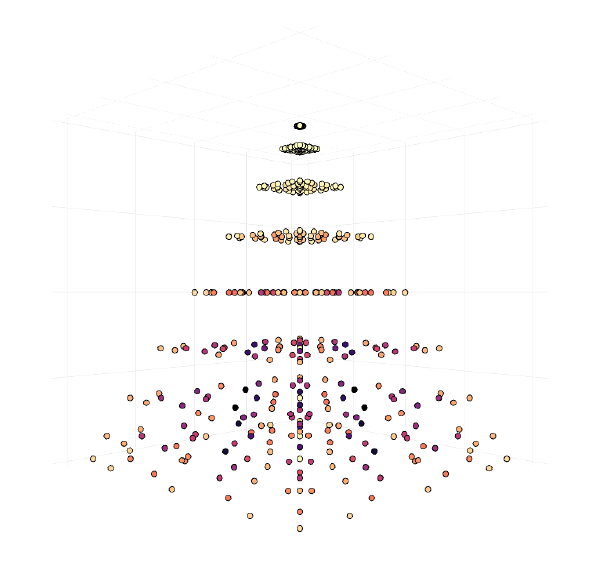}
        \caption{Algebraic}
        \label{fig:pyr-alg15}
      \end{subfigure}\hfill
      
      \caption{Geometric and algebraic initializations for degree 15 quadrature rule on the pyramid with darker nodes having larger weight values.}
      \label{fig:pyr-geoVSalg}
  \end{figure}

  Unlike the algebraic construction, the geometric one does not automatically satisfy the moment equations. Optimization is required to get a quadrature rule using the geometric construction, but a relatively low number of LMA iterations suffices. Nevertheless, it yields significantly fewer and more evenly distributed nodes, as can be seen in \cref{fig:pyr-geoVSalg}. Furthermore, the gap in node count between the two constructions increases with the quadrature degree. For degree 30, the geometric construction has 572 fewer nodes than the algebraic one.

  \section{Node Elimination}
  \label{node-elimination}
  Once an initial quadrature is constructed, the next step is to reduce the number of nodes as much as possible. Our approach combines orbit elimination and orbit collapse. The elimination of an orbit is the removal of the nodes associated with it. A critical aspect of this part is the order of elimination attempts, as it affects the number of nodes that can be removed. Our strategy is inspired by the one described in \cite{Worku2026} for the triangle and tetrahedron, which removes orbits strictly from larger symmetries first, then within each symmetry, removes orbits with smallest weights first. We also remove orbits in this order, but symmetries sharing the same number of parameters are bundled and relatively prioritized by assigning them priority numbers. Consider an orbit of Cartesian weight $w_C$ and exponential weight $w_e$ that is part of a $d$-parameter symmetry with priority number $p \geq 1$. The orbit's priority relative to other orbits in a $d$-parameter symmetry is set to $p_C=p\times w_C$ under Cartesian parameterization and $p_e=p' + w_e$ under exponential parameterization, where $p' = \log(p)/s$. This setup ensures consistent prioritization under parameterization changes. Indeed, $p_C=e^{sp_e}$ because $w_C=e^{sw_e}$ via \cref{eq:exp2}, so $p_{C1} < p_{C2} \iff p_{e1} < p_{e2}$. Each domain has a pair of symmetries of the same size and the same number of parameters (see \cref{tab:symmetries}), and we experimented with assigning those symmetries the same relative priority. However, prioritizing $S_2$ over $S_3$ on the square and pyramid, $S_5$ over $S_6$ on the cube, and $S_5$ over $S_4$ on the prism leads to the elimination of more nodes than otherwise. \cref{tab:symmetry-bundles} shows which symmetries are bundled (relatively prioritized) and the relative priorities within each bundle, and \cref{alg:orbit-elimination} outlines how orbit elimination works on a bundle of symmetries.

  \begin{table}[!t]
  \caption{Symmetry bundles \(\mathcal{SB}\) and corresponding priority bundles \(\mathcal{PB}\).}
      \label{tab:symmetry-bundles}
      \centering
      \begin{tabular}{@{} c | l @{}}
          square \quad &
          $\begin{aligned}
              \mathcal{SB} &= [[S_4],\ [S_3, S_2],\ [S_1]] \\
              \mathcal{PB} &= [[1],\ [10^5, 1],\ [1]]
          \end{aligned}$ \\[1.0em] 
          \hline\addlinespace
          cube \quad &
          $\begin{aligned}
              \mathcal{SB} &= [[S_7],\ [S_6, S_5],\ [S_4, S_3, S_2],\ [S_1]] \\
              \mathcal{PB} &= [[1],\ [10^5, 1],\ [1, 10^5, 10^{10}],\ [1]]
          \end{aligned}$ \\[1.0em]
          \midrule\addlinespace
          prism \quad &
          $\begin{aligned}
              \mathcal{SB} &= [[S_6],\ [S_5, S_4],\ [S_3, S_2],\ [S_1]] \\
              \mathcal{PB} &= [[1],\ [1, 10^5],\ [1, 10^5],\ [1]]
          \end{aligned}$ \\[1.0em]
          \midrule\addlinespace
          pyramid \quad &
          $\begin{aligned}
              \mathcal{SB} &= [[S_4],\ [S_3, S_2],\ [S_1]] \\
              \mathcal{PB} &= [[1],\ [10^5, 1],\ [1]]
          \end{aligned}$ \\
      \end{tabular}
  \end{table}

  \begin{algorithm}[!t]
      \caption{Procedure for eliminating orbits of a quadrature $\mathcal{Q}$ with degree $q$ from a set of symmetries $\mathcal{S}$ with corresponding priority numbers in $\mathcal{P}$.}
      \label{alg:orbit-elimination}
      \begin{algorithmic}[1]
          \Procedure{RemoveOrbits}{$\mathcal{Q}, q, \mathcal{S}, \mathcal{P}$}
              \State{OrbitPriorities $\gets$ $\Call{ComputeOrbitPriorities}{\mathcal{Q},q, \mathcal{S}, \mathcal{P}}$}
              \State{$i \gets 0$}
              \Comment{number of explored orbits}
              \While{$i <$ \#orbits from $\mathcal{S}$}
                  \If{$\mathcal{Q}$ has a single orbit}
                      \State{\textbf{break}}
                  \EndIf
                  
                  \State{TargetOrbit $\gets$ orbit with $i$th smallest priority}
                  \State{$\mathcal{Q}' \gets$ a copy of $\mathcal{Q}$ with Target removed}
                  \State{Residual $\gets \Call{LMA}{\mathcal{Q}',q}$ }
      
                  \If{Residual $<$ tolerance}
                      \State{$\mathcal{Q} \gets \mathcal{Q}'$}
                      \State{OrbitPriorities $\gets$ $\Call{ComputeOrbitPriorities}{\mathcal{Q}, \mathcal{S}, \mathcal{P}}$}
                      \State{$i \gets 0$}
                  \Else 
                      \If{Target's symmetry has been fully explored}
                      \State{remove it from $\mathcal{S}$ and its priority number from $\mathcal{P}$}
                      \EndIf
                      \State{$i \gets i+1$}
                      \Comment{moving onto the next orbit}
                  \EndIf
              
              \EndWhile
          \EndProcedure
      \end{algorithmic}
  \end{algorithm}

  The collapse of an orbit is its replacement by another orbit in a symmetry with fewer parameters. The sizes of symmetries decrease with the number of parameters, so a collapse reduces the number of nodes in a quadrature rule. Consider an orbit $O_i$ of weight $w$ in a $d$-parameter symmetry $S_i$. We attempt to collapse the orbit onto $(d-1)$-parameter symmetries within a prescribed geometric distance, which is referred to as the collapse threshold and denoted by $\mathcal{CT}$. Each symmetry of a domain corresponds to a region in it, hence the use of the term distance. Suppose $O_i$ is eligible for a collapse on the symmetry $S_f$. The replacement orbit $O_f$ is set to be the orbit in $S_f$ closest to $O_i$ with weight $(|S_i|/|S_f|)\times w$, where $|S_i|$ and $|S_f|$ are the sizes of $S_i$ and $S_f$, respectively. The total weight of the nodes associated with $O_i$ is the same as the total weight of the ones associated with $O_f$. In case $O_f$ can be collapsed on multiple symmetries, we first attempt the collapse on symmetries with lower priority numbers, targeting smallest size symmetries to maximize node elimination. The collapse threshold is set to $\mathcal{CT}=0.25$ for degrees below 31 on the square and below 20 on 3D domains where nodes are sparse. It is set to $0.1$ for higher degrees, where nodes are denser. ~\cref{alg:node-elimination} describes how orbit elimination and collapse are combined.

  \begin{algorithm}[!t]
      \caption{Procedure for eliminating nodes of a degree $q$ quadrature $\mathcal{Q}$ using orbit elimination and collapse given symmetry bundles $\mathcal{SB}$ and priority bundles $\mathcal{PB}$, as defined in Table~\ref{tab:symmetry-bundles}, and a collapse threshold $\mathcal{CT}$.} 
      \label{alg:node-elimination}
      \begin{algorithmic}[1]
          \Procedure{RemoveNodes}{$\mathcal{Q}, q, \mathcal{SB}, \mathcal{PB}, \mathcal{CT}$}
              \For{each index $i$ of $\mathcal{SB}$}
                  \Comment{indexing starts at 1}
                  \State{$\Call{RemoveOrbits}{\mathcal{Q},q,\mathcal{SB}[i],\mathcal{PB}[i]}$}
                  \Comment{see Algorithm~\ref{alg:orbit-elimination}}
                  \For{$j = 1:i$}
                      \State{$\Call{CollapseOrbits}{\mathcal{Q},q,\mathcal{SB}[j]}$}
                  \EndFor
              \EndFor
          \EndProcedure
      \end{algorithmic}
  \end{algorithm}

  Orbit collapse makes the node elimination process more robust by allowing a flow of orbits between symmetries. Suppose that the minimal node quadrature rule of degree $q$ for some domain has 5 orbits on its $S_k$ symmetry, but the initial quadrature rule only has 2 orbits on $S_k$. In this case, the initialization must have extra orbits on other symmetries since it is not minimal. Orbit elimination alone cannot lead to the minimal number of nodes, but with orbit collapse, it is possible to get there because $S_k$ can take orbits from other symmetries.

  An important part of the node elimination process is the stopping criteria of the LMA iteration. It depends on two parameters: the convergence check interval, $C_I$, and the maximum iteration, $M_I$. After $C_I\times k$ iterations, $k \in \mathbb{N}$, the algorithm checks whether the residual is less than the initial residual divided by $10^k$ and stops if that is not the case. This setup helps avoid wasting hundreds of iterations on orbits that cannot be eliminated. The algorithm also stops after $M_I$ iterations. $C_I$ is set to values between 20 and 70, increasing it with the number of variables, and $M_I$ is set to 500.

  \section{Results}
  \label{results}
  \subsection{Quadrature Rules}
  Our algorithm has led to numerous f-SPI quadrature rules with fewer nodes than existing ones across all four domains. Tables~\ref{tab:sqr-results},~\ref{tab:hex-results},~\ref{tab:pri-results}, and~\ref{tab:pyr-results} contain the number of nodes for each of our quadrature rules, with degrees where they are better than existing ones underlined. Existing quadrature rules have fewer nodes than presented in this work in only two instances (for degrees 7 and 9 rules on the prism), and the difference is exactly one node for each case. We only report odd degree rules on the square and cube because, similar to Legendre-Gauss rules in 1D, there are no even-degree fully symmetric quadrature rules on these domains. The nodal distributions for selected quadrature rules on the square, cube, prism, and pyramid are displayed in \cref{fig:sqr-results,fig:hex-results,fig:pri-results,fig:pyr-results}, respectively, in \ref{app1}.

  \begin{table}[!t]
      \centering
      \small
      \caption{Results for the square.}
      \label{tab:sqr-results}
      \begin{tabular}{c|cc||c|cc||c|c||c|c||c|c}
          $q$ & New & \cite{WITHERDEN20151232} & $q$ & New & \cite{WITHERDEN20151232} & $q$ & New & $q$ & New & $q$ & New \\ \hline
          1 & 1 & 1 & \underline{17} & 57 & 60 & \underline{33} & 201 & \underline{49} & 433 & \underline{65} & 757 \\
          3 & 4 & 4 & 19 & 72 & 72 & \underline{35} & 224 & \underline{51} & 465 & \underline{67} & 805 \\
          5 & 8 & 8 & 21 & 85 & 85 & \underline{37} & 249 & \underline{53} & 501 & \underline{69} & 849 \\
          7 & 12 & 12 & \underline{23} & 101 &  & \underline{39} & 277 & \underline{55} & 541 & \underline{71} & 904  \\
          9 & 20 & 20 & \underline{25} & 120 &  & \underline{41} & 304 & \underline{57} & 576 & \underline{73} & 953 \\
          11 & 28 & 28 & \underline{27} & 137 &  & \underline{43} & 336 & \underline{59} & 613 & \underline{75} & 1001 \\
          13 & 37 & 37 & \underline{29} & 157 &  & \underline{45} & 365 & \underline{61} & 660 & \underline{77} & 1049 \\
          15 & 48 & 48 & \underline{31} & 177 &  & \underline{47} & 397 & \underline{63} & 709 &  &  
      \end{tabular}
  \end{table}

  \begin{table}[!t]
      \centering
      \small
      \caption{Results for the cube.}
      \label{tab:hex-results}
      \begin{tabular}{c|cc||c|cc||c|cc||c|c||c|c}
          $q$ & New & \cite{jaskowiec2021addendum} & $q$ & New & \cite{jaskowiec2021addendum} & $q$ & New & \cite{jaskowiec2021addendum} & $q$ & New & $q$ & New  \\ \hline
          1 & 1 & 1 \cite{WITHERDEN20151232} & 11 & 90 & 90 & \underline{21} & 495 & 580 & \underline{31} & 1478 & \underline{41} & 3338 \\
          3 & 8 &  & \underline{13} & 148 & 154 & \underline{23} & 617 &  & \underline{33} & 1787 & \underline{43} & 3870 \\
          5 & 14 & 14 & \underline{15} & 199 & 256 & \underline{25} & 828 &  & \underline{35} & 2102 & \underline{45} &  4414 \\
          7 & 34 & 34 & \underline{17} & 282 & 346 & \underline{27} & 984 &  & \underline{37} & 2506 &  &  \\
          9 & 58 & 58 & \underline{19} & 369 & 454 & \underline{29} & 1258 &  & \underline{39} & 2856 &  &  
      \end{tabular}
  \end{table}

  \begin{table}[!t]
      \centering
      \small
      \caption{Results for the prism.}
      \label{tab:pri-results}
      \begin{tabular}{c|cc||c|cc||c|cc||c|c||c|c}
          $q$ & New & \cite{jaskowiec2021addendum} & $q$ & New & \cite{jaskowiec2021addendum} & $q$ & New & \cite{jaskowiec2021addendum} & $q$ & New & $q$ & New  \\ \hline
          1&  1&  1 \cite{WITHERDEN20151232}&  8&  46&  46& \underline{15}& 220& 238& \underline{22}& 635& \underline{29}& 1311\\
          2&  5&  5&  9&  60&  59& \underline{16}& 264&  287& \underline{23}& 695& \underline{30}& 1460\\
          3&  8&  8&  \underline{10}&  82&  84& \underline{17}& 299&  338& \underline{24}& 772& & \\
          4&  11&  11&  \underline{11}&  97&  99& \underline{18}& 359&  396& \underline{25}& 878&  &   \\
          5&  16&  16&  \underline{12}&  127&  136& \underline{19}& 405&  420& \underline{26}& 986&  &  \\
          6& 28& 28& \underline{13}& 146& 162& \underline{20}& 470& 518& \underline{27}& 1080& &\\
          7& 36& 35& \underline{14}& 182& 194& \underline{21}& 545& & \underline{28}& 1197& &\\\end{tabular}
  \end{table}

  \begin{table}
      \centering
      \small
      \caption{Results for the pyramid}
      \label{tab:pyr-results}
      \begin{tabular}{c|cc||c|cc||c|cc||c|c||c|c}
          $q$ & New & \cite{Jaskowiec2021main} & $q$ & New & \cite{Jaskowiec2021main} & $q$ & New & \cite{Jaskowiec2021main} & $q$ & New & $q$ & New  \\ \hline
          1&  1&  1 \cite{WITHERDEN20151232} &  \underline{8}&  44&  47& \underline{15}& 203& 234& \underline{22}& 595& \underline{29}& 1275\\
          2&  5&  5&  \underline{9}&  56&  62& \underline{16}& 254&  285& \underline{23}& 665& \underline{30}& 1417\\
          3&  6&  6&  \underline{10}&  76&  80& \underline{17}& 293&  319& \underline{24}& 756& & \\
          4&  10&  10&  \underline{11}&  92&  103& \underline{18}& 344&  357& \underline{25}& 843&  &   \\
          5&  15&  15&  \underline{12}&  120&  127& \underline{19}& 394&  418& \underline{26}& 942&  &  \\
          6& 23& 23& \underline{13}& 142& 152& \underline{20}& 464& 489& \underline{27}& 1050& &\\
          7& 31& 31& \underline{14}& 175& 184& \underline{21}& 522& & \underline{28}& 1171& &\\
      \end{tabular}
  \end{table}

  We quantify the efficiency of a quadrature rule for a given domain by comparing it with the best we can build using quadrature rules on lower-dimensional domains. In the case of the square, cube, and prism, the reference rules are tensor-product constructions. These coincide with our initial rules, except where we intentionally added extra orbits to make sure some symmetries are non-empty. For the pyramid, we use algebraic constructions as a reference without any extra orbits added. The efficiency of a degree $q$ quadrature rule with $n$ nodes, relative to a reference quadrature rule with $n_r$ nodes, is defined as
  \begin{equation}
      e = \frac{n_r-n}{n_r}.
  \end{equation}
  The number of reference nodes is given by 
  \begin{align*}
      &\left(\left\lceil \frac{q+1}{2}\right\rceil\right)^2, \ \left(\left\lceil\frac{q+1}{2}\right\rceil\right)^3, \\ &\left(\left\lceil\frac{q+1}{2}\right\rceil\right)\times n_{T}(q), \text{ and }  \left(\left\lceil\frac{(q+2)+1}{2}\right\rceil\right)\times n_{S}(q)
  \end{align*}
  on the square, cube, prism, and pyramid, respectively, where $n_{T}(q)$ and $n_S(q)$ are the number of nodes of the best degree $q$ quadrature rule on the triangle \cite{Worku2026} and square (our results), respectively.

  The significance of this metric lies not in its value at a particular degree but rather how it changes with the quadrature degree. For low degrees, the optimization problem has few variables and is easy to solve; we could conveniently explore various node reduction strategies and are thus confident that the quadrature rules obtained there have an almost minimal number of nodes. By low degrees, we mean 1-31 and 1-21 for the square and cube, respectively, and 1-14 for the prism and pyramid, because it takes a very short amount of time to go through the initialization and node elimination process. We use the efficiency trend at low degrees to make a judgment about the performance of the algorithm at high degrees. 

  \begin{figure}[!t]
      \begin{subfigure}{0.5\linewidth}
        \centering
        \includegraphics[
          width=\linewidth, 
        ]{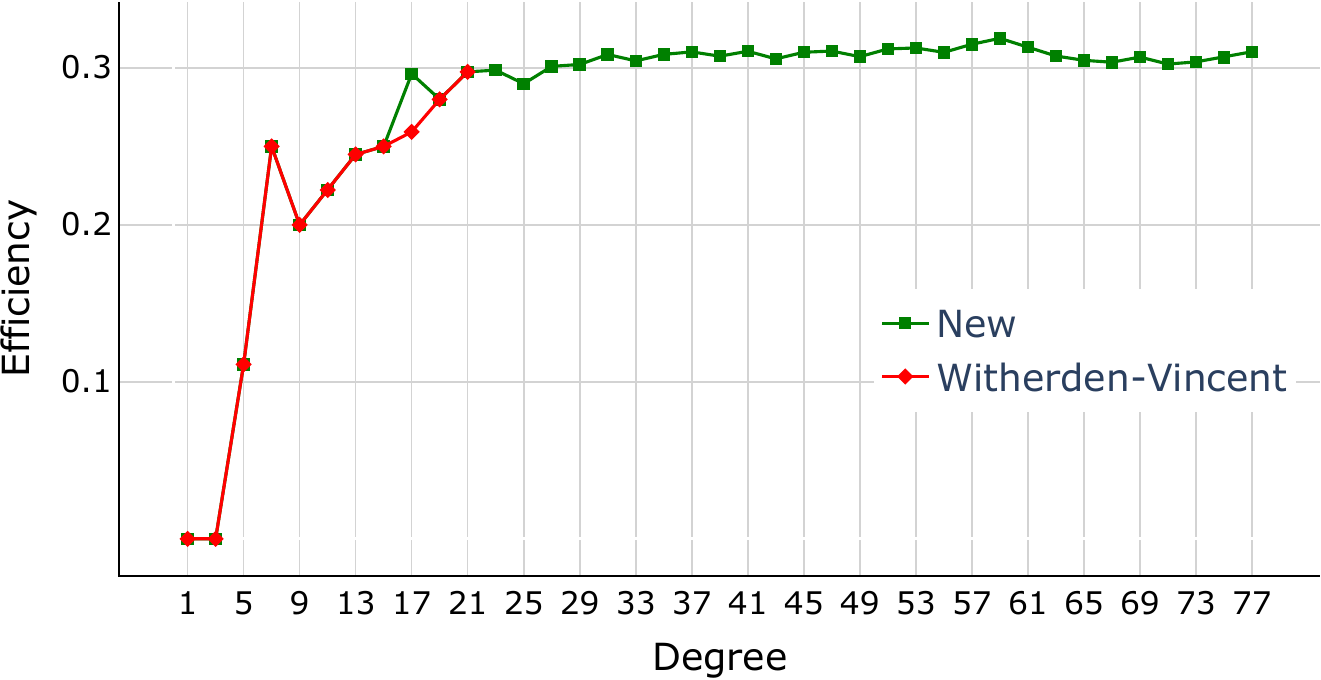}
        \caption{Square.}
        \label{fig:sqr-eff}
      \end{subfigure}\hfill
      \begin{subfigure}{0.5\linewidth}
        \centering
        \includegraphics[
          width=\linewidth, 
        ]{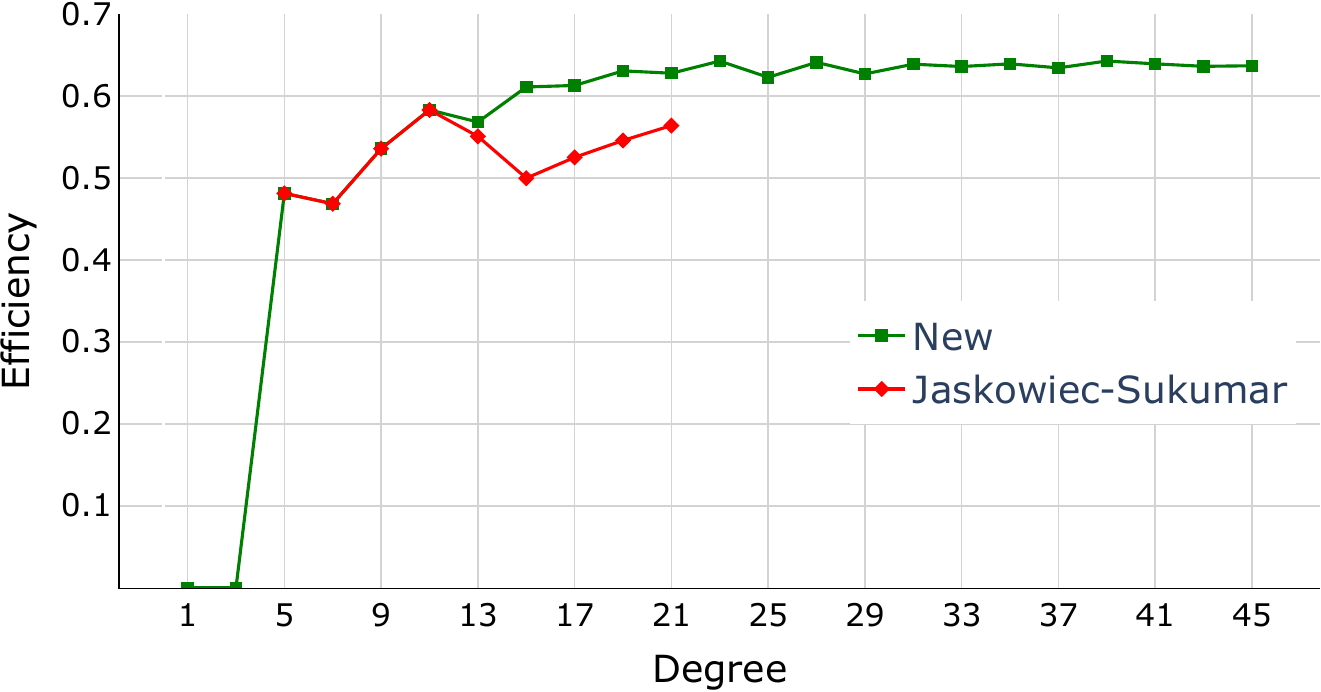}
        \caption{Cube.}
        \label{fig:hex-eff}
      \end{subfigure}
      
      \caption{Efficiency on the square and cube.}
      \label{fig:sqr-hex-eff}
  \end{figure}

  \begin{figure}[!t]
      \begin{subfigure}{0.5\linewidth}
        \centering
        \includegraphics[
          width=\linewidth, 
        ]{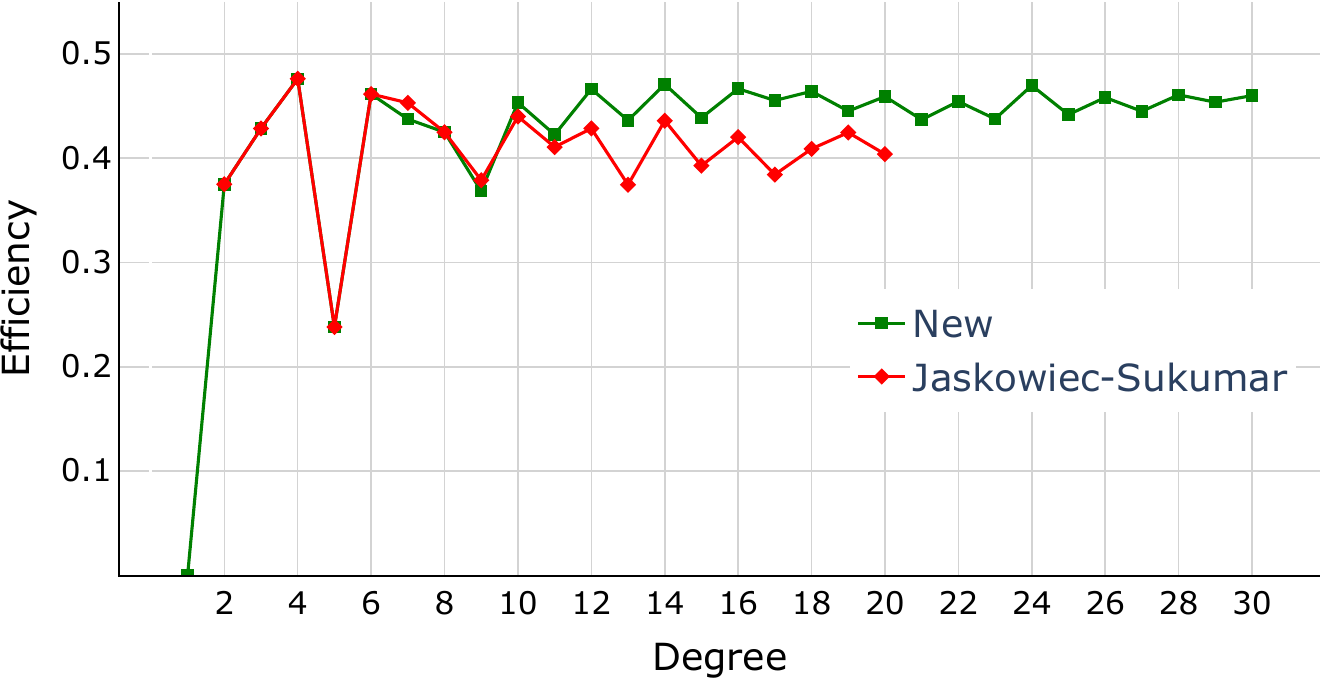}
        \caption{Prism.}
        \label{fig:pri-eff}
      \end{subfigure}\hfill
      \begin{subfigure}{0.5\linewidth}
        \centering
        \includegraphics[
          width=\linewidth, 
        ]{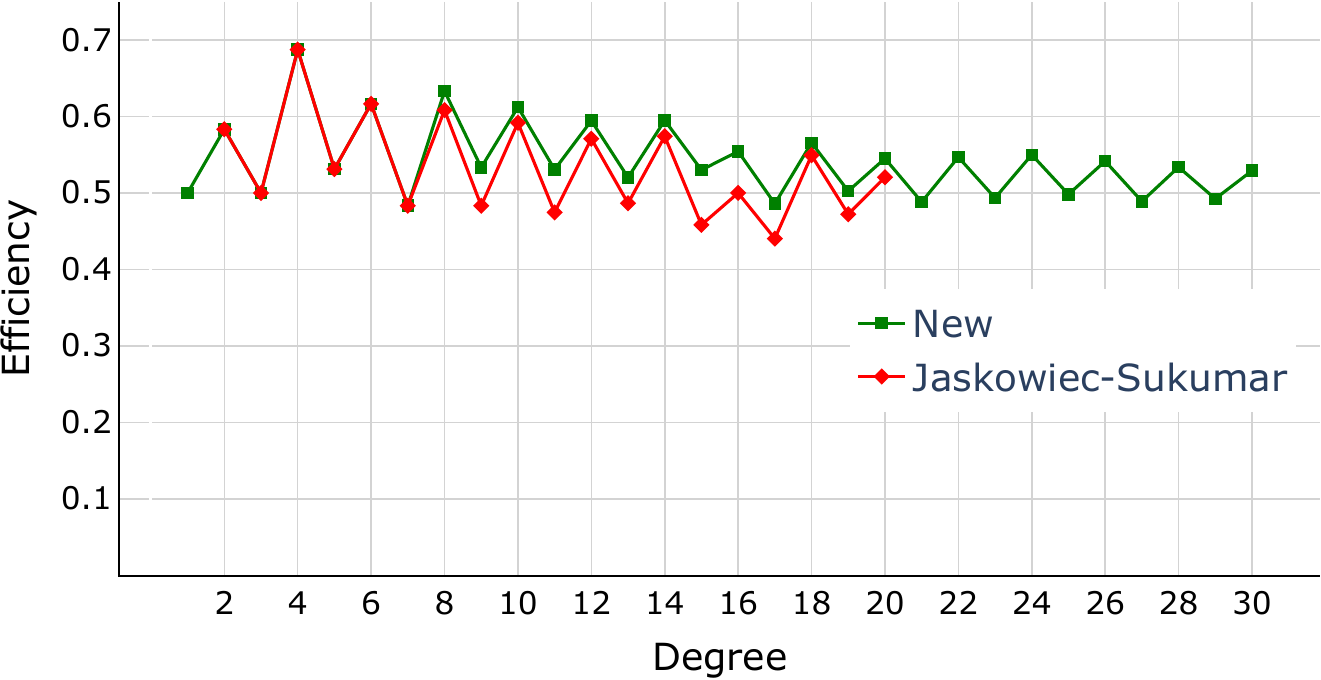}
        \caption{Pyramid.}
        \label{fig:pyr-eff}
      \end{subfigure}
      
      \caption{Efficiency on the prism and pyramid.}
      \label{fig:pri-pyr-eff}
  \end{figure}

  On the square, the efficiency increases rapidly between $q=1$ and $q=21$, going from 0 to about 0.30, but does not change much beyond that, see Figure~\ref{fig:sqr-eff}. The efficiency at degree 77 being comparable to that of degrees 21-31 is a good indicator of a steady performance of the node elimination algorithm. There is a similar trend on the cube (prism), see Figure~\ref{fig:hex-eff} (\ref{fig:pri-eff}), where the efficiency reaches 0.61 (0.45) at degree 15 (10) and stays within $0.625\pm 0.025$ ($0.44\pm 0.02$) all the way up to degree 45 (30). It also has a steady trend on the pyramid, see Figure~\ref{fig:pyr-eff}, but an oscillatory one, unlike the other three domains.

  \subsection{Numerical Experiments}
  The quadrature rules presented in this work are based on polynomials, but they can be used to approximate integrals of all kinds of functions. Oscillatory functions are usually a challenge, and we use them to test our results. Specifically, we use the functions
  \begin{align}
      \label{eq:test-functions}
      f(x,y) &= cos(k_1x+k_2y) \quad \text{and}
      \\ 
      g(x,y,z) &= cos(k_1x+k_2y+k_3z)
  \end{align}
  in two and three dimensions, respectively, where $k_1$, $k_2$, and $k_3$ are whole numbers. We proceed by approximating the integral
  \begin{equation}
      \label{eq:test-integral-2d}
      I(k_1,k_2)=\int_{-1}^1\int_{-1}^1 f(x,y)dxdy = \frac{4\sin(k_1)\sin(k_2)}{k_1k_2}
  \end{equation}
  using quadrature rules on a mesh of squares and approximating the integral
  \begin{equation}
      \label{eq:test-integral-3d}
      I(k_1,k_2,k_3)=\int_{-1}^1\int_{-1}^1\int_{-1}^1 g(x,y,z)dxdydz = \frac{8\sin(k_1)\sin(k_2)\sin(k_3)}{k_1k_2k_3}
  \end{equation}
  using quadrature rules on a mesh of cubes, prisms, and pyramids.

  \begin{figure}[!t]
        \centering
        \includegraphics[
          width=0.8\linewidth, 
        ]{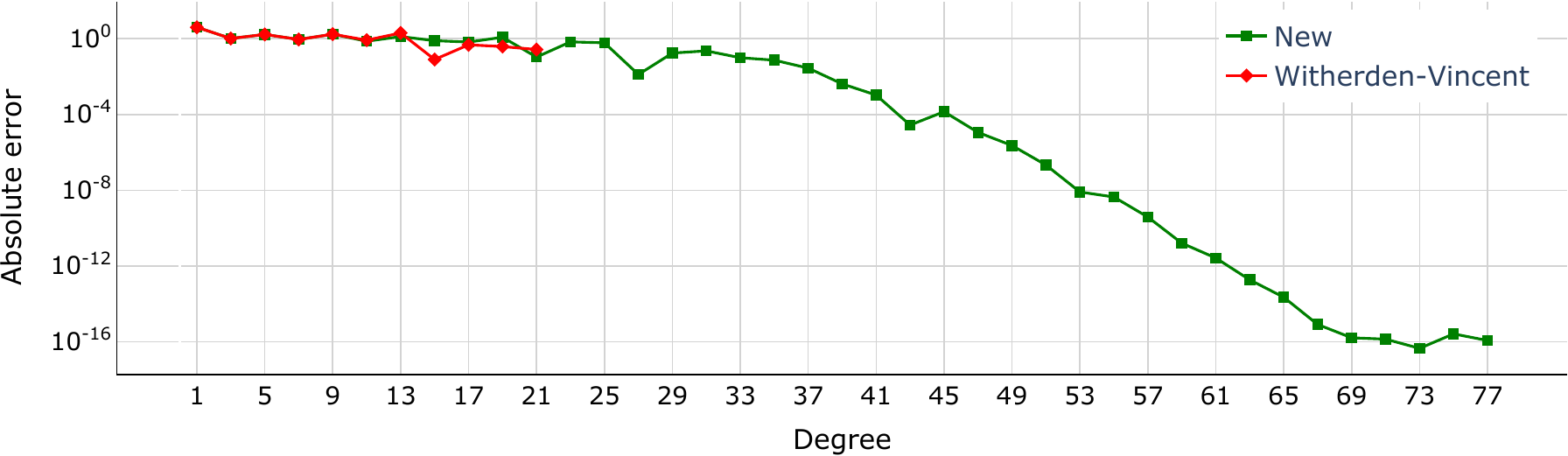}
        \caption{Approximation of $I(10,15)$ with mesh of one square.}
        \label{fig:sqr-degrees-a10b25}
  \end{figure}

  \begin{figure}[!t]
        \centering
        \includegraphics[
          width=0.8\linewidth, 
        ]{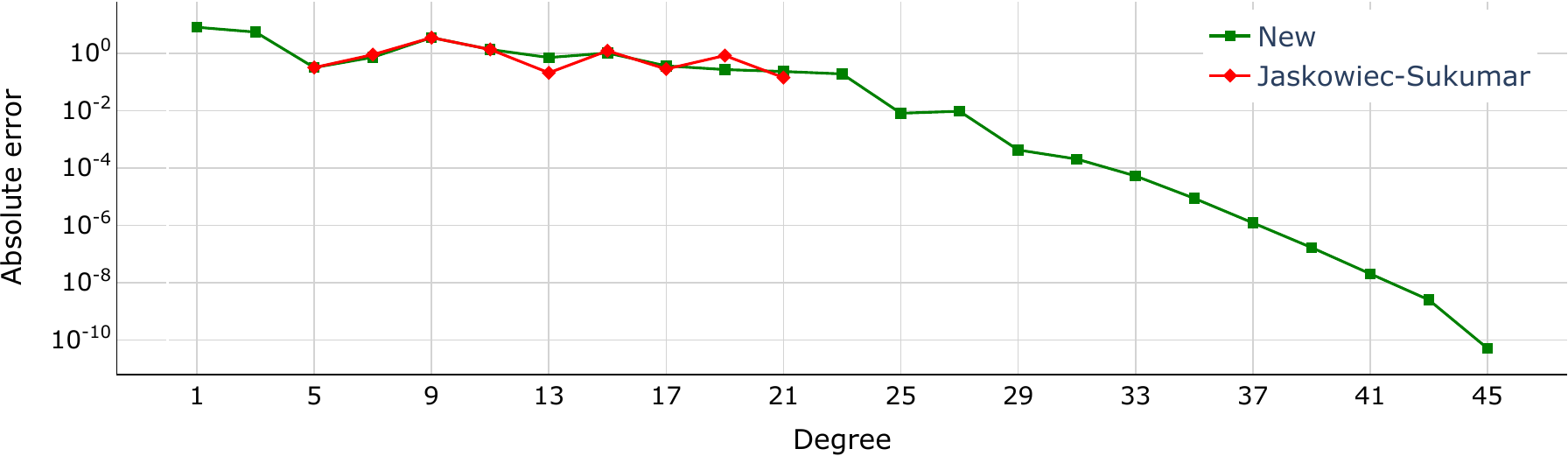}
        \caption{Approximation of $I(5,5,15)$ with mesh of one cube.}
        \label{fig:hex-degrees-a5b5c15}
  \end{figure}

  \begin{figure}[!t]
        \centering
        \includegraphics[
          width=0.8\linewidth, 
        ]{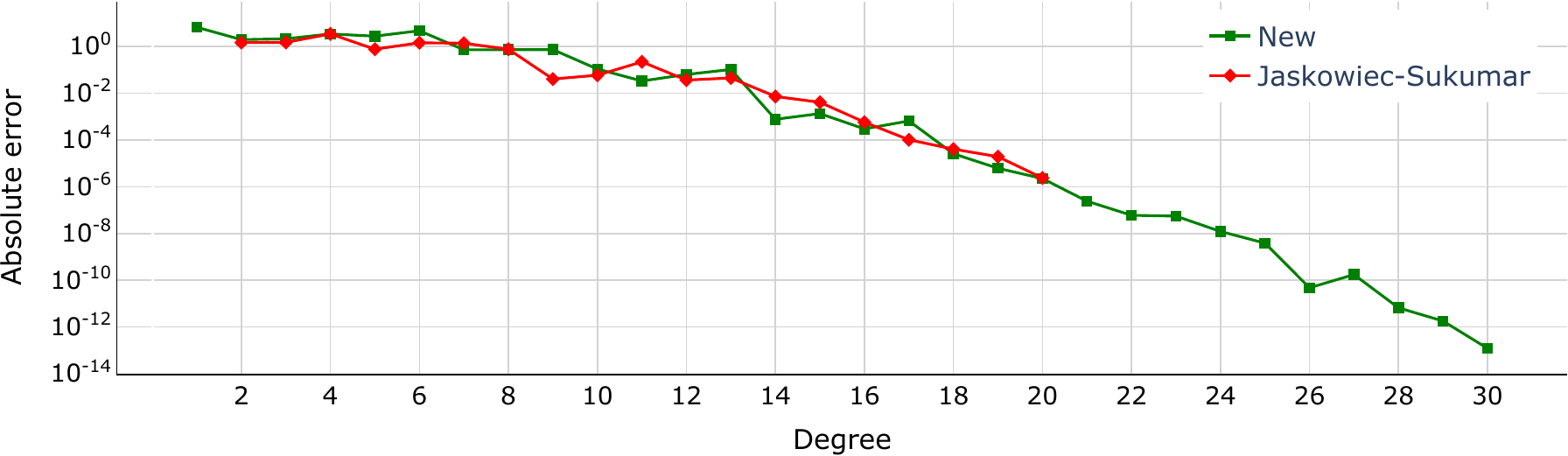}
        \caption{Approximation of $I(1,1,10)$ with mesh of two prisms.}
        \label{fig:pri-degrees-a1b1c10}
  \end{figure}

  \begin{figure}[!t]
        \centering
        \includegraphics[
          width=0.8\linewidth, 
        ]{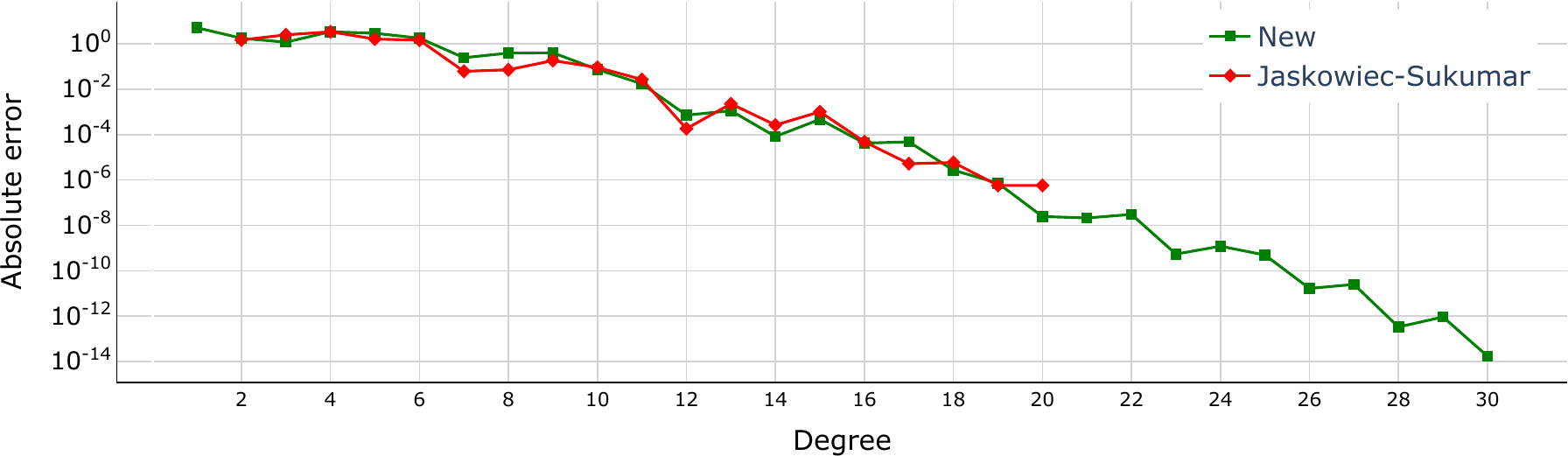}
        \caption{Approximation of $I(1,1,10)$ with mesh of six pyramids.}
        \label{fig:pyr-degrees-a1b1c10}
  \end{figure}

\begin{figure}[!t]
      \begin{subfigure}{0.475\linewidth}
        \centering
        \includegraphics[
          width=\linewidth, 
        ]{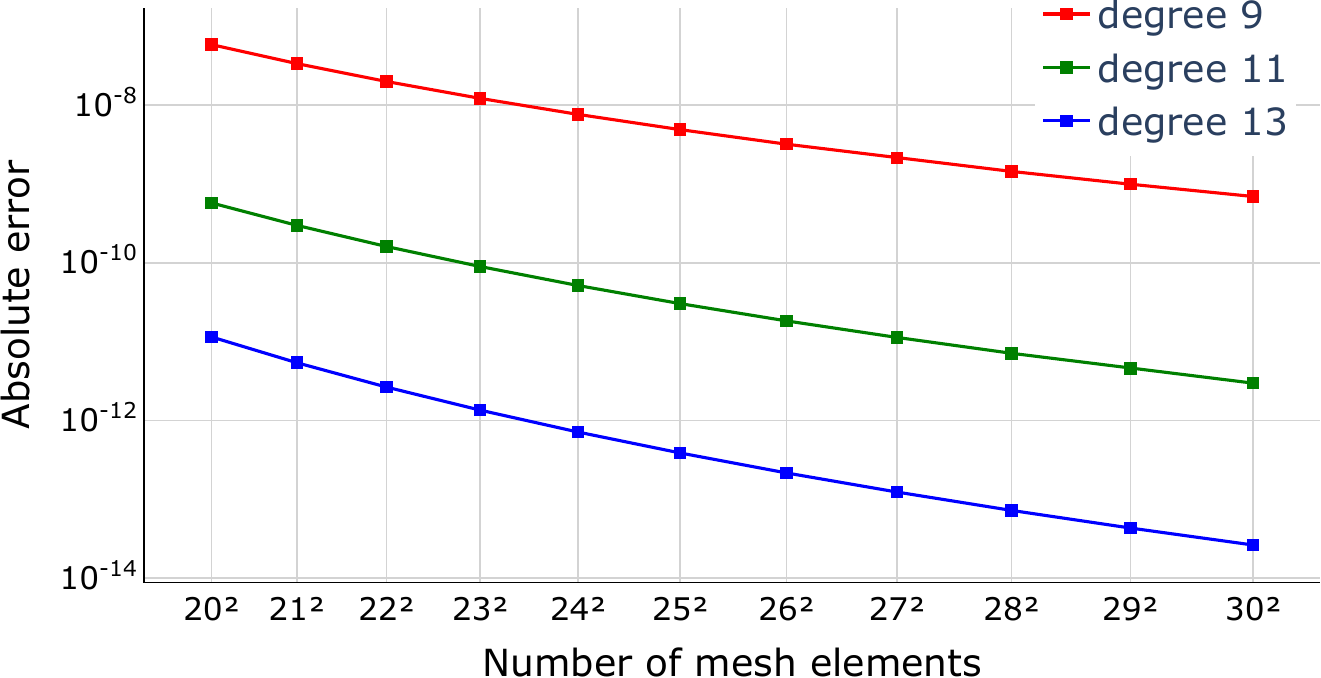}
        \caption{Absolute error.}
        \label{fig:sqr-mesh-error}
      \end{subfigure}\hfill
      \begin{subfigure}{0.475\linewidth}
        \centering
        \includegraphics[
          width=\linewidth, 
        ]{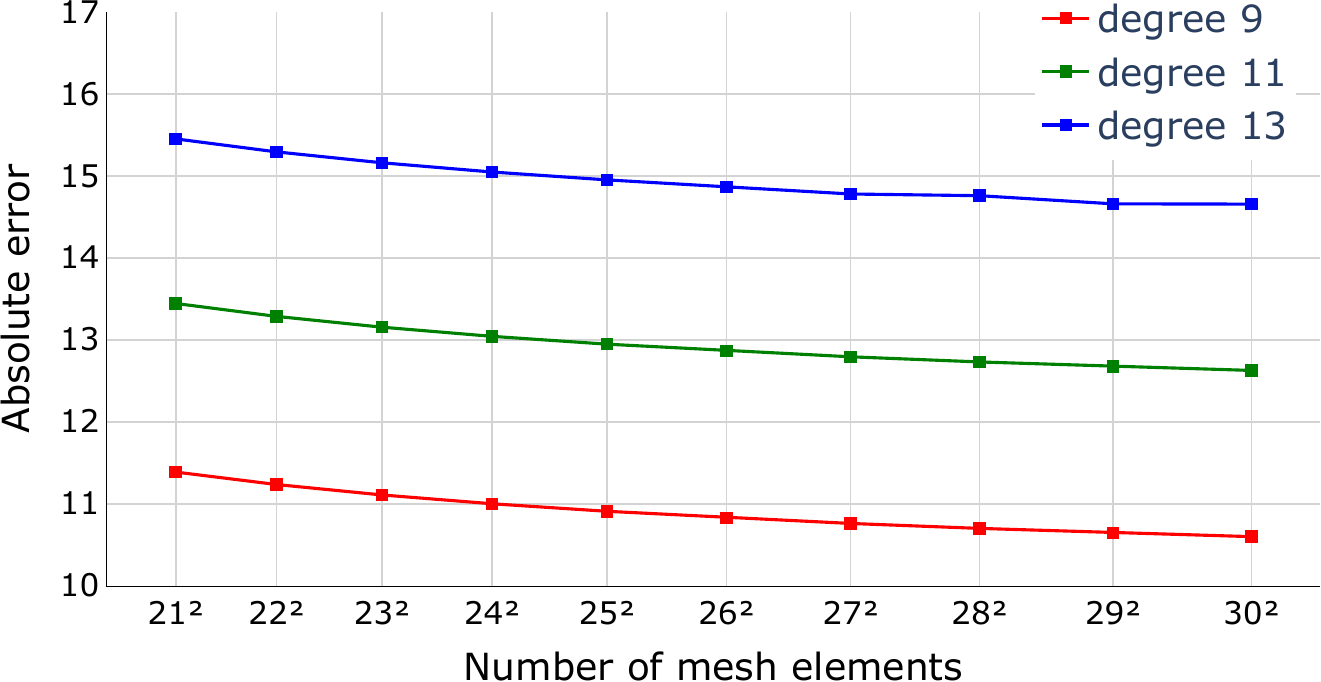}
        \caption{Convergence rate.}
        \label{fig:sqr-mesh-rate}
      \end{subfigure}
      
      \caption{Approximation of $I(30,30)$ with meshes of squares of various number of elements.}
      \label{fig:sqr-mesh-test}
  \end{figure}

\begin{figure}[!t]
      \begin{subfigure}{0.475\linewidth}
        \centering
        \includegraphics[
          width=\linewidth, 
        ]{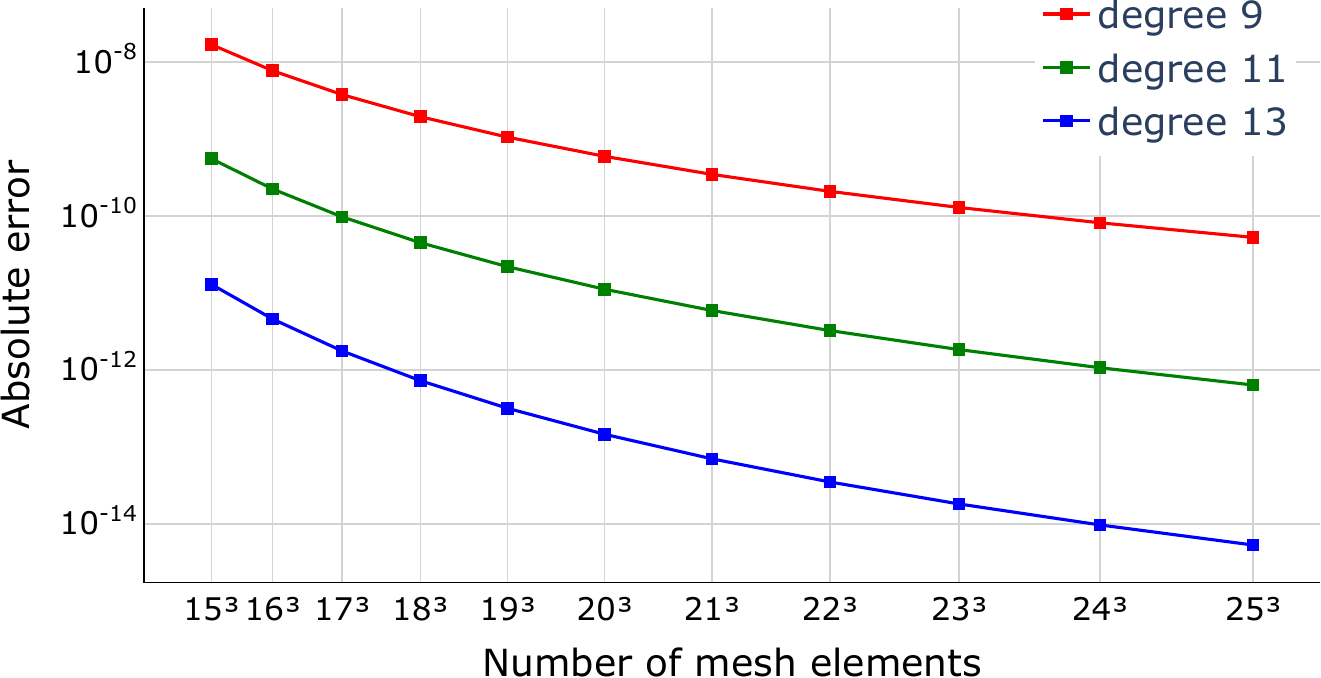}
        \caption{Absolute error.}
        \label{fig:hex-mesh-error}
      \end{subfigure}\hfill
      \begin{subfigure}{0.475\linewidth}
        \centering
        \includegraphics[
          width=\linewidth, 
        ]{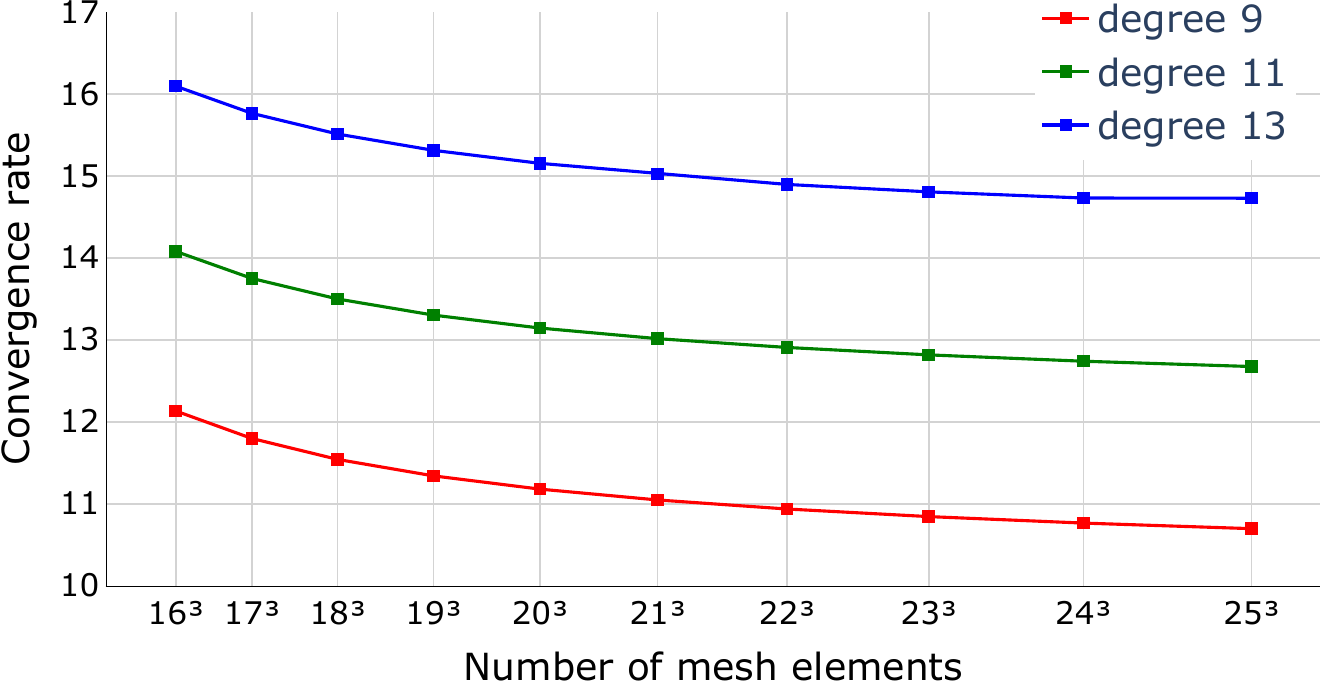}
        \caption{Convergence rate.}
        \label{fig:hex-mesh-rate}
      \end{subfigure}
      
      \caption{Approximation of $I(30,10,20)$ with meshes of cubes of various number of elements.}
      \label{fig:hex-mesh-test}
  \end{figure}

\begin{figure}[!t]
      \begin{subfigure}{0.475\linewidth}
        \centering
        \includegraphics[
          width=\linewidth, 
        ]{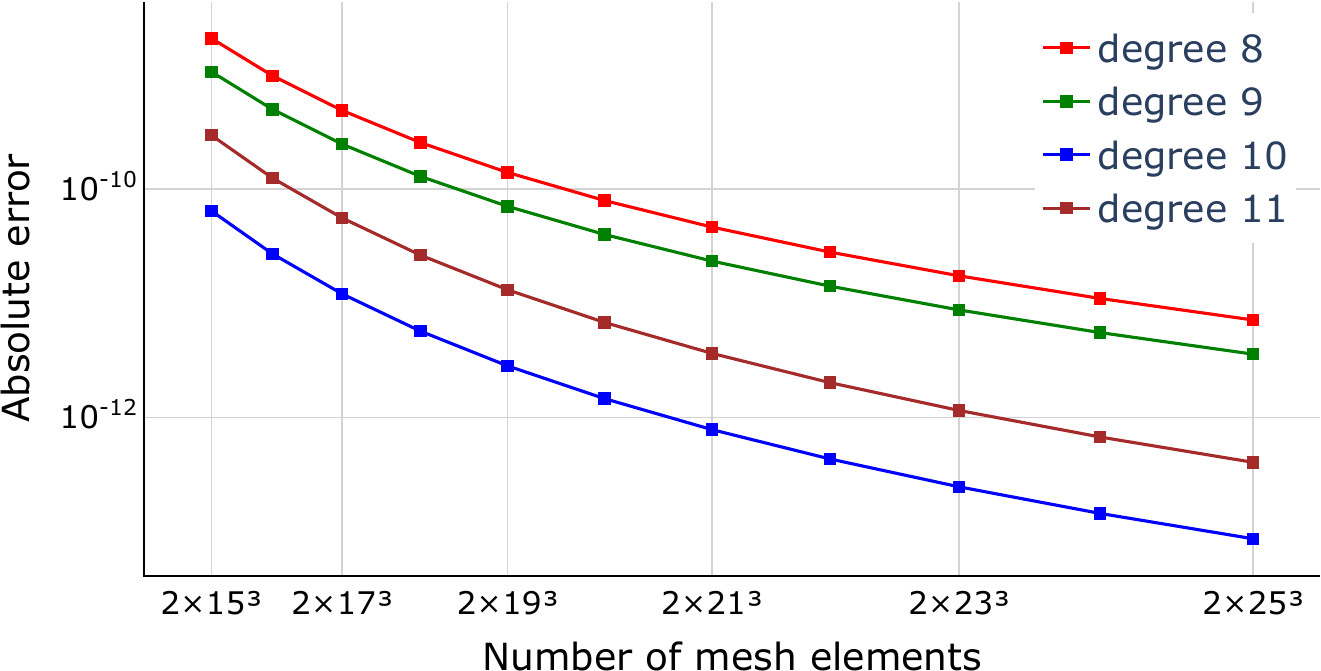}
        \caption{Absolute error.}
        \label{fig:pri-mesh-error}
      \end{subfigure}\hfill
      \begin{subfigure}{0.475\linewidth}
        \centering
        \includegraphics[
          width=\linewidth, 
        ]{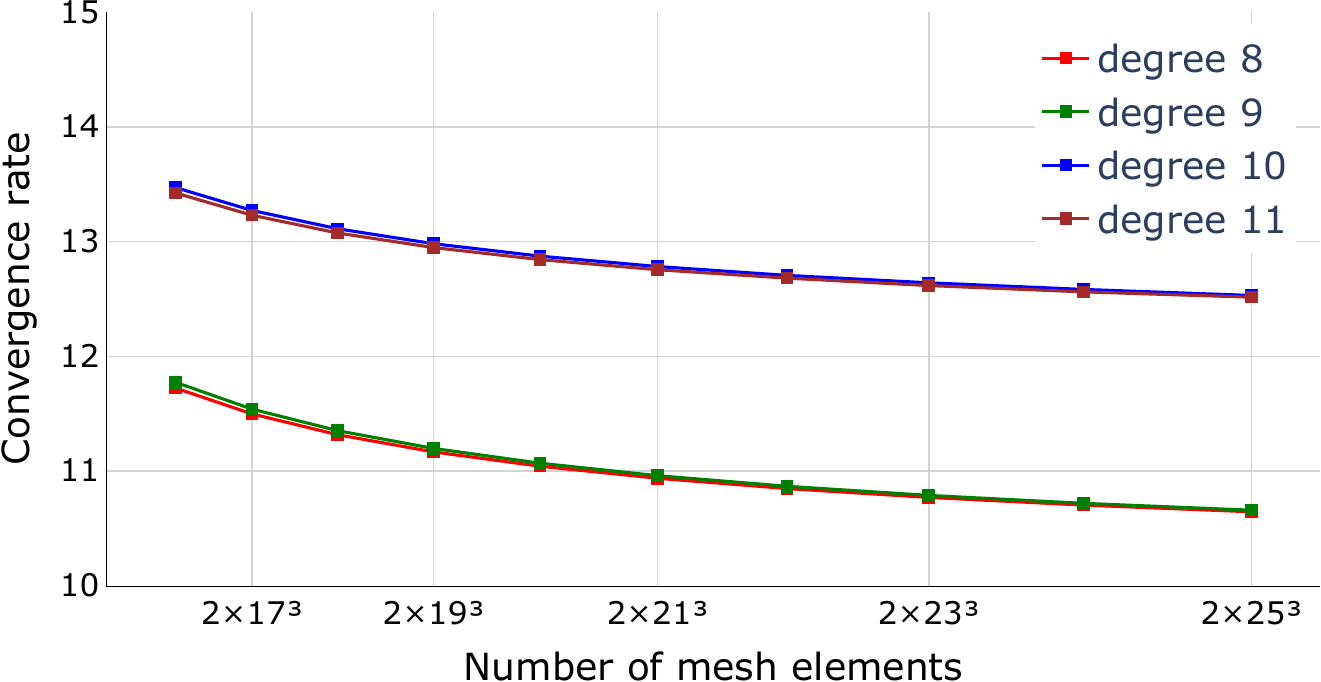}
        \caption{Convergence rate.}
        \label{fig:pri-mesh-rate}
      \end{subfigure}
      
      \caption{Approximation of $I(40,20,20)$ with meshes of prisms of various number of elements.}
      \label{fig:pri-mesh-test}
  \end{figure}

\begin{figure}[!t]
      \begin{subfigure}{0.475\linewidth}
        \centering
        \includegraphics[
          width=\linewidth, 
        ]{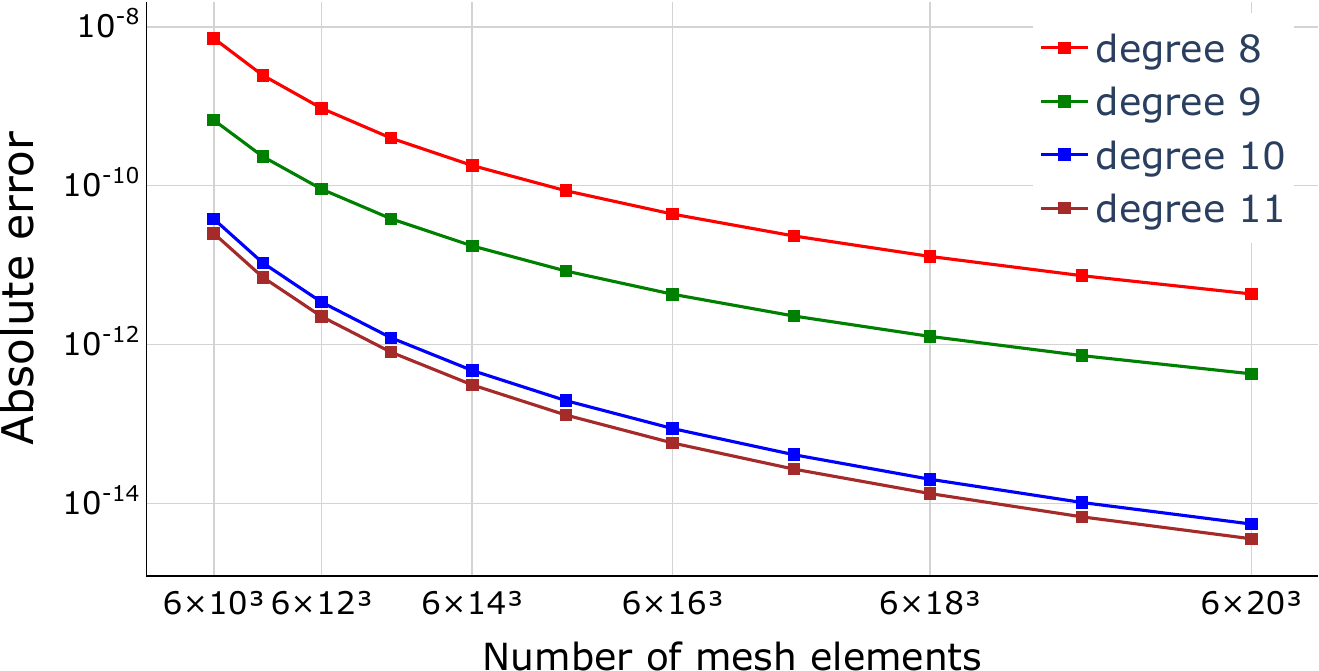}
        \caption{Absolute error.}
        \label{fig:pyr-mesh-error}
      \end{subfigure}\hfill
      \begin{subfigure}{0.475\linewidth}
        \centering
        \includegraphics[
          width=\linewidth, 
        ]{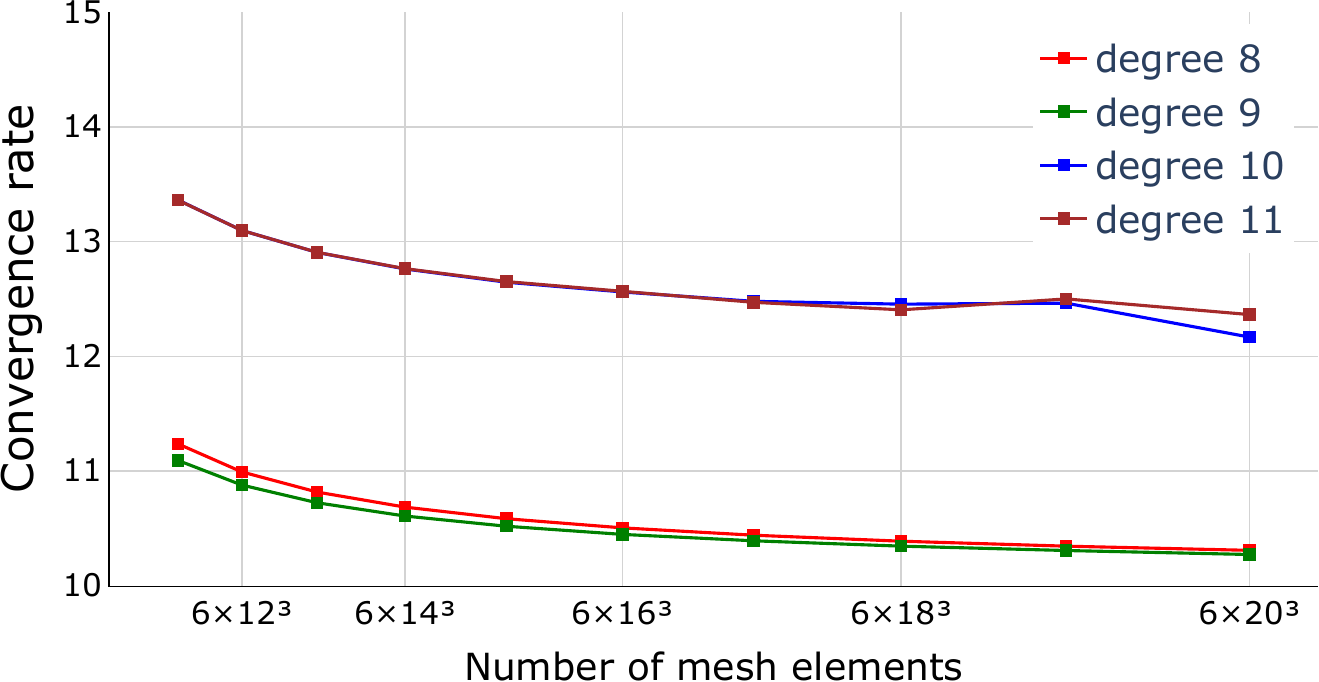}
        \caption{Convergence rate.}
        \label{fig:pyr-mesh-rate}
      \end{subfigure}
      
      \caption{Approximation of $I(15,5,15)$ with meshes of pyramids of various number of elements.}
      \label{fig:pyr-mesh-test}
  \end{figure}

  Figures~\ref{fig:sqr-degrees-a10b25},~\ref{fig:hex-degrees-a5b5c15},~\ref{fig:pri-degrees-a1b1c10}, and~\ref{fig:pyr-degrees-a1b1c10} show that the absolute errors in the approximations of $I(k_1,k_2)$ and $I(k_1,k_2,k_3)$ using our quadrature rules decrease as the quadrature degree increases --- at rates comparable to those of existing quadrature rules. In addition, Figures~\ref{fig:sqr-mesh-test},~\ref{fig:hex-mesh-test},~\ref{fig:pri-mesh-test}, and~\ref{fig:pyr-mesh-test} display a convergence rate of at least $q+1$ for a degree $q$ quadrature rule with mesh refinement. The convergence rate is close to $q+2$ for even degrees, and similar observations were made in \cite{Worku2026} for the tetrahedra.

  \section{Conclusions}
  \label{conclusions}

  In this work, we developed a practical framework for constructing f‑SPI quadrature rules on commonly used multidimensional reference elements. The approach combines a variable parameterization that enforces positivity and interiority with a Levenberg–Marquardt optimization algorithm and a symmetry‑aware node‑reduction strategy based on orbit elimination and orbit collapse. Orbit collapse, in particular, improves robustness by allowing orbits to transition between symmetry types, enabling node reduction beyond what orbit elimination alone typically permits. Using these ingredients, we produced new f‑SPI quadrature rules on the square, cube, prism, and pyramid. The resulting rules reach degrees up to 77 on the square, 45 on the cube, and 30 on the prism and pyramid, and for most degrees we obtain improved node counts relative to previously published f‑SPI rules. Numerical experiments confirm the expected convergence behavior and demonstrate accuracy comparable to existing rules. Complete node and weight data are provided to facilitate direct use in high‑order discretizations of PDEs. Several directions for future work remain open. Further improvements in initialization strategies and solver robustness could extend the attainable degrees in three dimensions and reduce node counts further. Extending the approach to construct quadrature rules with nodes on element boundaries is also of interest, as such rules enable efficient entropy‑stable discretizations within the summation‑by‑parts framework.

  \clearpage
  \appendix
  \section{Nodal Distributions of Selected Quadrature Rules}
  \label{app1}

  \begin{figure}[!htbp]
      \begin{subfigure}{0.5\linewidth}
        \centering
        \includegraphics[
          width=\linewidth, 
        ]{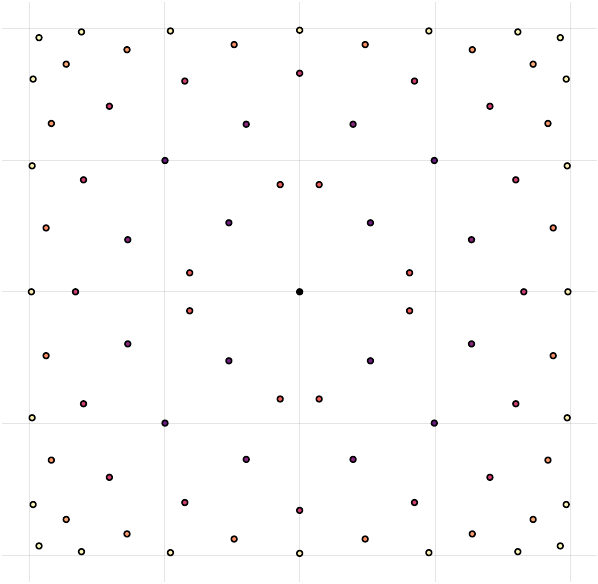}
        \caption{Square: degree 21}
        \label{fig:sqr-q21}
      \end{subfigure}\hfill
      \begin{subfigure}{0.5\linewidth}
        \centering
        \includegraphics[
          width=\linewidth, 
        ]{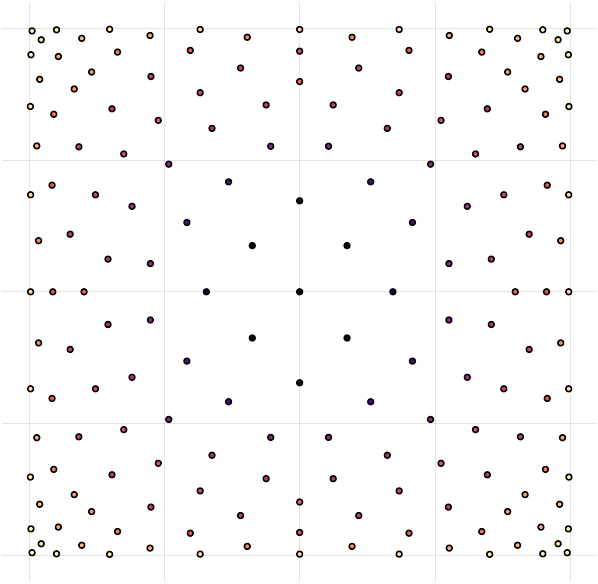}
        \caption{Square: degree 31}
        \label{fig:sqr-q31}
      \end{subfigure}

      \begin{subfigure}{0.5\linewidth}
        \centering
        \includegraphics[
          width=\linewidth, 
        ]{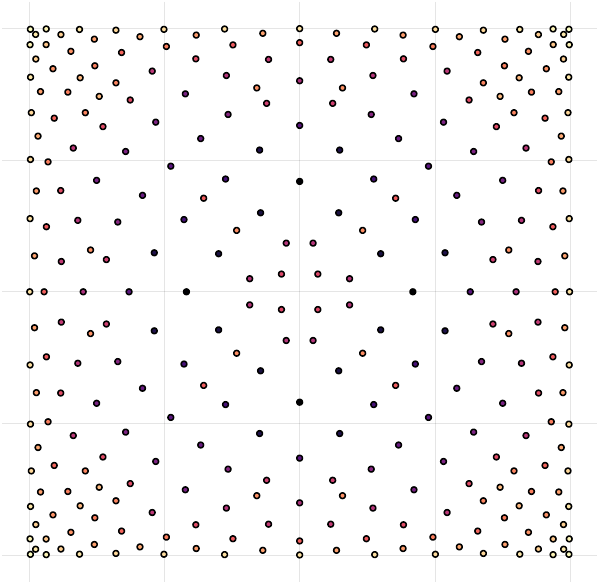}
        \caption{Square: degree 41}
        \label{fig:sqr-q41}
      \end{subfigure}\hfill
      \begin{subfigure}{0.5\linewidth}
        \centering
        \includegraphics[
          width=\linewidth, 
        ]{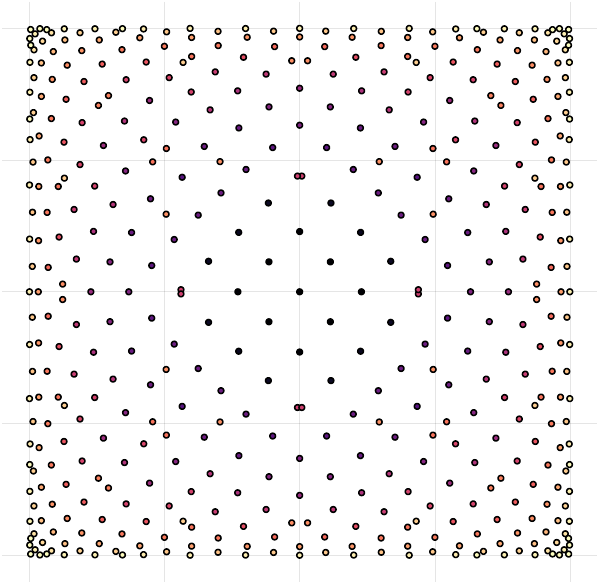}
        \caption{Square: degree 51}
        \label{fig:sqr-q51}
      \end{subfigure}
      
      \caption{Final quadrature rules obtained on the square for degrees 21, 31, 41, and 51; darker nodes have more weight.}
      \label{fig:sqr-results}
  \end{figure}

  \begin{figure}[!thbp]
      \begin{subfigure}{0.5\linewidth}
        \centering
        \includegraphics[
          width=\linewidth, 
        ]{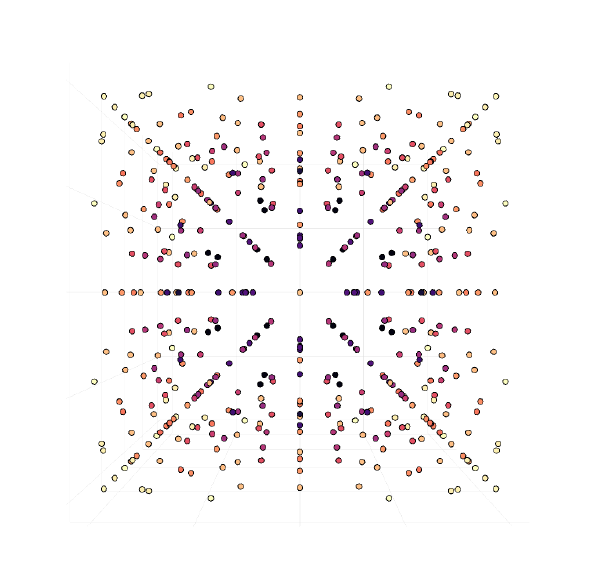}
        \caption{Cube: degree 21}
        \label{fig:hex-q21-Fview}
      \end{subfigure}\hfill
      \begin{subfigure}{0.5\linewidth}
        \centering
        \includegraphics[
          width=\linewidth, 
        ]{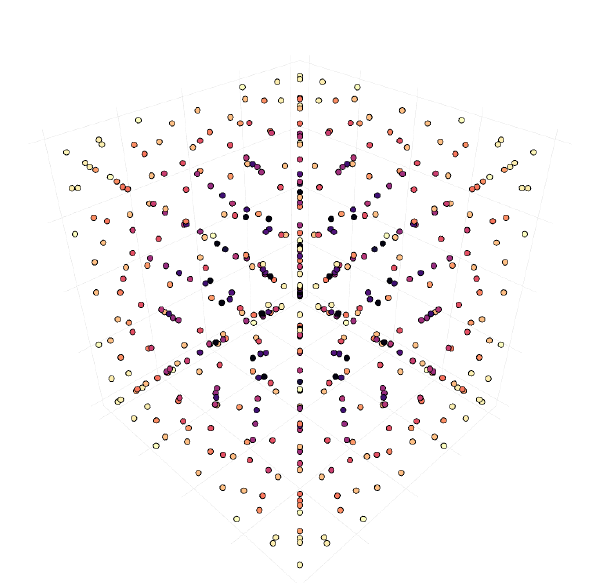}
        \caption{Cube: degree 21}
        \label{fig:hex-q21-Dview}
      \end{subfigure}

      \begin{subfigure}{0.5\linewidth}
        \centering
        \includegraphics[
          width=\linewidth, 
        ]{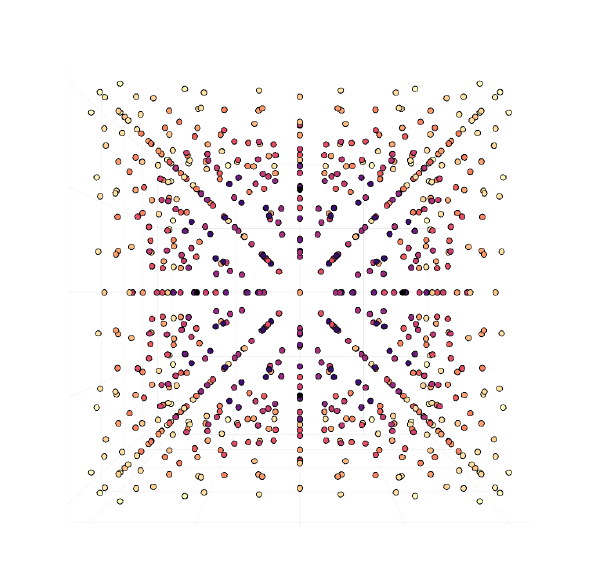}
        \caption{Cube: degree 25}
        \label{fig:hex-q25-Fview}
      \end{subfigure}\hfill
      \begin{subfigure}{0.5\linewidth}
        \centering
        \includegraphics[
          width=\linewidth, 
        ]{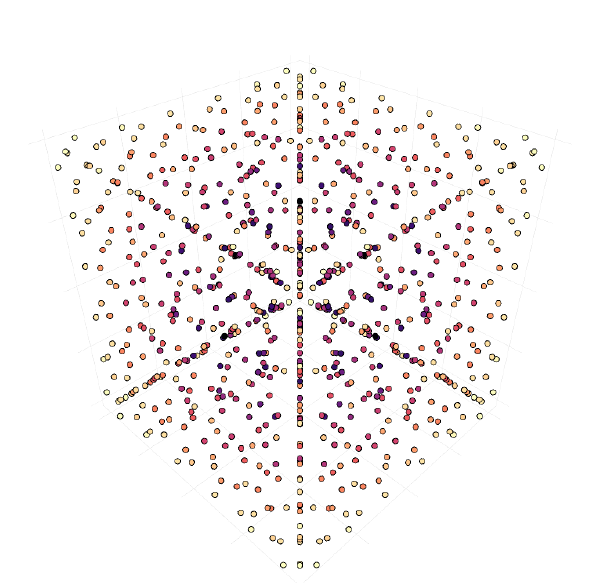}
        \caption{Cube: degree 25}
        \label{fig:hex-q25-Dview}
      \end{subfigure}
      
      \caption{Final quadrature rules obtained on the cube for degrees 21 and 25, with a face view (left) and a diagonal view (right); darker nodes have more weight.}
      \label{fig:hex-results}
  \end{figure}

  \begin{figure}[!thbp]
      \centering
      \begin{subfigure}{0.5\linewidth}
        \centering
        \includegraphics[
          width=\linewidth, 
        ]{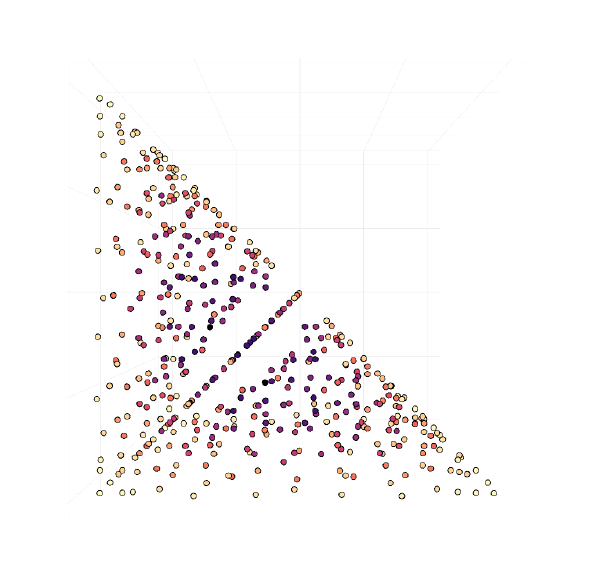}
        \caption{Prism: degree 20}
        \label{fig:hex-q20-Tview}
      \end{subfigure}\hfill
      \begin{subfigure}{0.5\linewidth}
        \centering
        \includegraphics[
          width=\linewidth, 
        ]{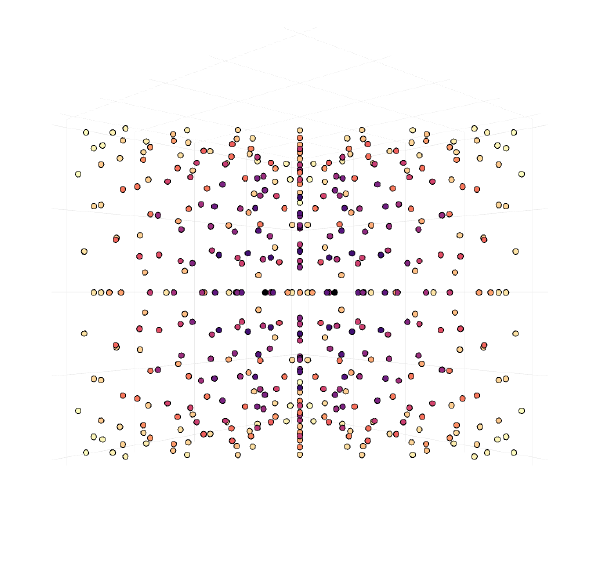}
        \caption{Prism: degree 20}
        \label{fig:pri-q20-Dview}
      \end{subfigure}

      \begin{subfigure}{0.5\linewidth}
        \centering
        \includegraphics[
          width=\linewidth, 
        ]{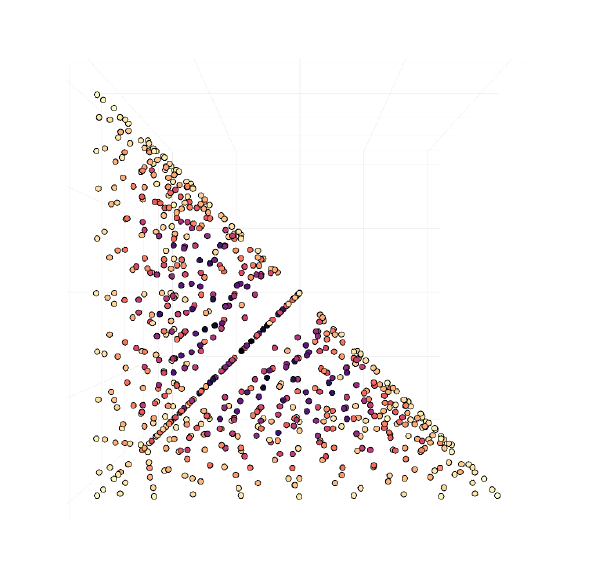}
        \caption{Prism: degree 23}
        \label{fig:pri-q23-Tview}
      \end{subfigure}\hfill
      \begin{subfigure}{0.5\linewidth}
        \centering
        \includegraphics[
          width=\linewidth, 
        ]{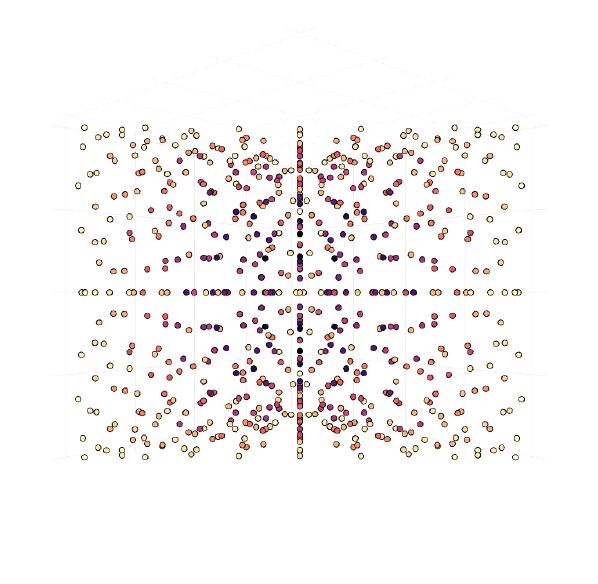}
        \caption{Prism: degree 23}
        \label{fig:pri-q23-Dview}
      \end{subfigure}
      
      \caption{Final quadrature rules obtained on the prism for degrees 20 and 23 with a top view (left) and a diagonal view (right); darker nodes have more weight.}
      \label{fig:pri-results}
  \end{figure}

  \begin{figure}[!thbp]
      \centering
      \begin{subfigure}{0.5\linewidth}
        \centering
        \includegraphics[
          width=\linewidth, 
        ]{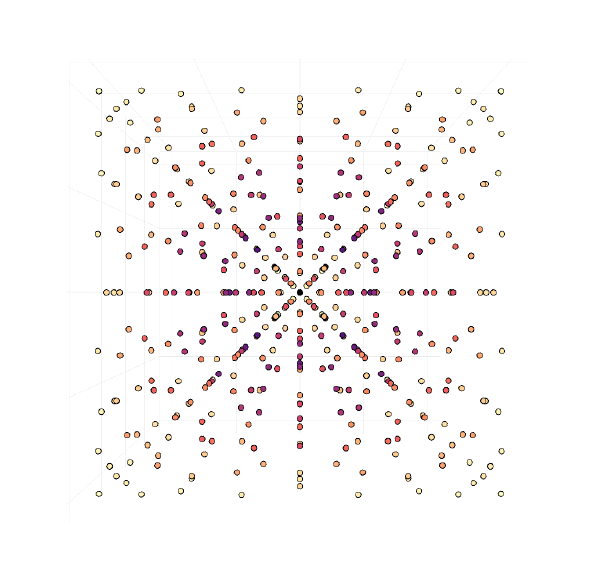}
        \caption{Pyramid: degree 20}
        \label{fig:pyr-q20-Bview}
      \end{subfigure}\hfill
      \begin{subfigure}{0.5\linewidth}
        \centering
        \includegraphics[
          width=\linewidth, 
        ]{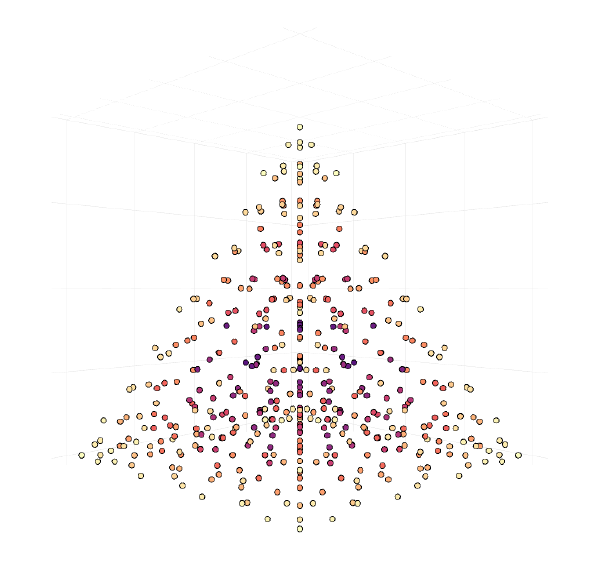}
        \caption{Pyramid: degree 20}
        \label{fig:pyr-q20-xView}
      \end{subfigure}

      \begin{subfigure}{0.5\linewidth}
        \centering
        \includegraphics[
          width=\linewidth, 
        ]{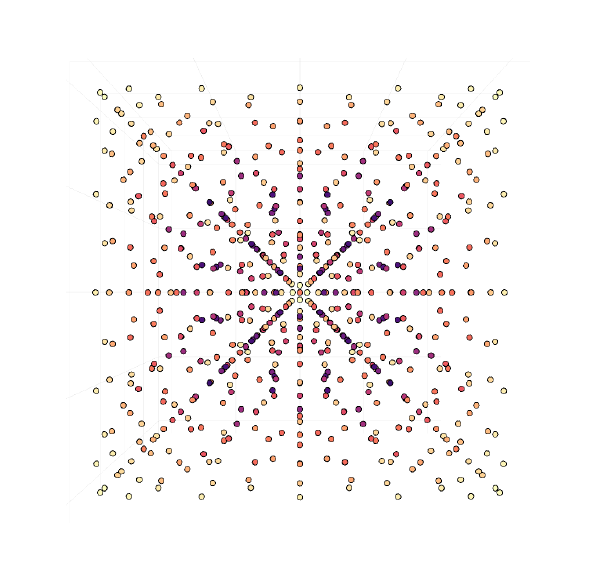}
        \caption{Pyramid: degree 23}
        \label{fig:pyr-q23-Bview}
      \end{subfigure}\hfill
      \begin{subfigure}{0.5\linewidth}
        \centering
        \includegraphics[
          width=\linewidth,
        ]{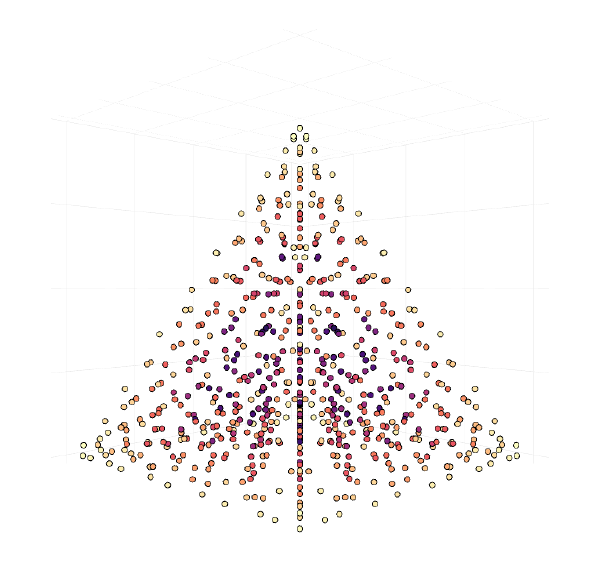}
        \caption{Pyramid: degree 23}
        \label{fig:pyr-q23-xView}
      \end{subfigure}
      
      \caption{Final quadrature rules obtained on the pyramid for degrees 20 and 23 with a base view (left) and a side view (right); darker nodes have more weight.}
      \label{fig:pyr-results}
  \end{figure}

  \section{Supplementary Materials}\label{app2}
  The nodal and weight data of the quadrature rules derived in this work can be found at \url{https://github.com/mdiallo-fula/SymmetricPositiveInteriorCubatures.jl/tree/master/data/published-results}.

  \bibliographystyle{elsarticle-num}
  \bibliography{references}

@conference{osti_4491151,
      author       = {Reed, W H and Hill, T R},
      title        = {Triangular mesh methods for the neutron transport equation},
      place        = {United States},
      organization = {Los Alamos Scientific Lab., N.Mex. (USA)},
      year         = {1973},
      month        = {10}
}

@article{COCKBURN1998199,
    title = {The {R}unge–{K}utta Discontinuous {G}alerkin Method for Conservation Laws V: Multidimensional Systems},
    journal = {Journal of Computational Physics},
    volume = {141},
    number = {2},
    pages = {199-224},
    year = {1998},
    issn = {0021-9991},
    author = {Bernardo Cockburn and Chi-Wang Shu}
}

@book{hesthaven2008,
    author={Hesthaven, Jan S and Warburton, Tim},
    title={Nodal discontinuous {G}alerkin methods: algorithms, analysis, and applications},
    year={2008},
    publisher={Springer}
}

@article{KIRBY2003249,
    title = {De-aliasing on non-uniform grids: algorithms and applications},
    journal = {Journal of Computational Physics},
    volume = {191},
    number = {1},
    pages = {249-264},
    year = {2003},
    issn = {0021-9991},
    author = {Robert M. Kirby and George Em Karniadakis},
}

@article{PERSSON20091585,
    title = {Discontinuous {G}alerkin solution of the {N}avier–{S}tokes equations on deformable domains},
    journal = {Computer Methods in Applied Mechanics and Engineering},
    volume = {198},
    number = {17},
    pages = {1585-1595},
    year = {2009},
    issn = {0045-7825},
    author = {P.-O. Persson and J. Bonet and J. Peraire},
}

@article{Williams2019EntropyL2DG,
    author  = {Williams, David M.},
    title   = {An analysis of discontinuous {G}alerkin methods for the compressible {E}uler equations: entropy and {$L_2$} stability},
    journal = {Numerische Mathematik},
    year    = {2019},
    volume  = {141},
    number  = {4},
    pages   = {1079--1120},
}

@article{chuluunbaatar2022,
    author  = {Chuluunbaatar, G. and Chuluunbaatar, O. and Gusev, A. A. and Vinitsky, S. I.},
    title   = {{PI}-type fully symmetric quadrature rules on the 3-, 6-simplexes},
    journal = {Computers and Mathematics with Applications},
    volume  = {124},
    pages   = {89--97},
    year    = {2022}
}

@article{FESTA20124296,
    title = {Computing almost minimal formulas on the square},
    journal = {Journal of Computational and Applied Mathematics},
    volume = {236},
    number = {17},
    pages = {4296-4302},
    year = {2012},
    issn = {0377-0427},
    author = {Mattia Festa and Alvise Sommariva},
}

@article{jaskowiec2021addendum,
    title={Addendum to the paper “{H}igh-order symmetric cubature rules for tetrahedra and pyramids”},
    author={Ja{\'s}kowiec, Jan and Sukumar, N},
    journal={International Journal for Numerical Methods in Engineering},
    volume={122},
    number={7},
    pages={1875--1883},
    year={2021},
    publisher={Wiley Online Library}
}

@article{Jaskowiec2021main,
    author = {Jaśkowiec, Jan and Sukumar, N.},
    title = {High-order symmetric cubature rules for tetrahedra and pyramids},
    journal = {International Journal for Numerical Methods in Engineering},
    volume = {122},
    number = {1},
    pages = {148-171},
    year = {2021}
}

@article{kubatko2013,
  author  = {Kubatko, Ethan J. and Yeager, Benjamin A. and Maggi, Ashley L.},
  title   = {New computationally efficient quadrature formulas for triangular prism elements},
  journal = {Computers \& Fluids},
  volume  = {73},
  pages   = {187--201},
  year    = {2013}
}

@article{SLOBODKINS2023229,
title = {A node elimination algorithm for cubature of high-dimensional polytopes},
journal = {Computers \& Mathematics with Applications},
volume = {144},
pages = {229-236},
year = {2023},
issn = {0898-1221},
author = {Arkadijs Slobodkins and Johannes Tausch},
}

@article{WANDZURAT20031829,
    title = {Symmetric quadrature rules on a triangle},
    journal = {Computers \& Mathematics with Applications},
    volume = {45},
    number = {12},
    pages = {1829-1840},
    year = {2003},
    issn = {0898-1221},
    author = {S. Wandzurat and H. Xiao},
}

@article{WITHERDEN20151232,
title = {On the identification of symmetric quadrature rules for finite element methods},
journal = {Computers \& Mathematics with Applications},
volume = {69},
number = {10},
pages = {1232-1241},
year = {2015},
issn = {0898-1221},
author = {F.D. Witherden and P.E. Vincent},
}

@article{Worku2026,
    author = {Worku, Zelalem Arega and Hicken, Jason E. and Zingg, David W.},
    title = {Very high-order symmetric positive-interior quadrature rules on triangles and tetrahedra},
    year = {2026},
    issue_date = {Jan 2026},
    publisher = {Elsevier Science Publishers B. V.},
    address = {NLD},
    volume = {472},
    number = {C},
    issn = {0377-0427},
    journal = {J. Comput. Appl. Math.}
}

@article{XIAO2010663,
    title = {A numerical algorithm for the construction of efficient quadrature rules in two and higher dimensions},
    journal = {Computers \& Mathematics with Applications},
    volume = {59},
    number = {2},
    pages = {663-676},
    year = {2010},
    issn = {0898-1221},
    author = {Hong Xiao and Zydrunas Gimbutas},
}

@article{duffy1982,
    author = {Duffy,  Michael G.},
    title = {Quadrature Over a Pyramid or Cube of Integrands with a Singularity at a Vertex},
    journal = {SIAM Journal on Numerical Analysis},
    volume = {19},
    number = {6},
    pages = {1260-1262},
    year = {1982},
}

@article{fernandez2014review,
    title     = {Review of summation-by-parts operators with simultaneous approximation terms for the numerical solution of partial differential equations},
    author    = {Del Rey Fern{\'a}ndez, David C  and Hicken, Jason E and Zingg, David W},
    journal   = {Computers \& Fluids},
    volume    = {95},
    pages     = {171--196},
    year      = {2014},
    publisher = {Elsevier}
}

@article{svard2014review,
    title     = {Review of summation-by-parts schemes for initial--boundary-value problems},
    author    = {Sv{\"a}rd, Magnus and Nordstr{\"o}m, Jan},
    journal   = {Journal of Computational Physics},
    volume    = {268},
    pages     = {17--38},
    year      = {2014},
    publisher = {Elsevier}
}

@article{hicken2016multidimensional,
    title     = {Multidimensional summation-by-parts operators: General theory and application to simplex elements},
    author    = {Hicken, Jason E and Del Rey Fern{\'a}ndez, David C and Zingg, David W},
    journal   = {SIAM Journal on Scientific Computing},
    volume    = {38},
    number    = {4},
    pages     = {A1935--A1958},
    year      = {2016},
    publisher = {SIAM}
}

@article{worku2024quadrature,
    title     = {Quadrature rules on triangles and tetrahedra for multidimensional summation-by-parts operators},
    author    = {Worku, Zelalem Arega and Hicken, Jason E and Zingg, David W},
    journal   = {Journal of Scientific Computing},
    volume    = {101},
    number    = {1},
    pages     = {24},
    year      = {2024},
    publisher = {Springer}
}
  \addcontentsline{toc}{section}{\refname}

\end{document}